\documentclass[reqno,12pt]{amsart}
\usepackage{latexsym}
\usepackage{amsxtra}
\usepackage{verbatim}
\usepackage{color}
\usepackage{amsthm,amsmath,amsfonts,amssymb,amscd,mathrsfs,graphics}
\usepackage{paralist}
\usepackage{amssymb}
\usepackage{times,float}
\usepackage{appendix}
\usepackage{txfonts}
\usepackage{hyperref,supertabular}
\usepackage{ifpdf}
\usepackage{enumerate}
\usepackage{fancyhdr}

\usepackage[top=1.5in, bottom=1.5in, left=1.5in, right=1.5in]{geometry}

\allowdisplaybreaks

\newtheorem{theorem}{Theorem}[section]
\newtheorem{lemma}[theorem]{Lemma}

\newtheorem{corollary}[theorem]{Corollary}
\newtheorem{definition}[theorem]{Definition}
\newtheorem{example}[theorem]{Example}
\newtheorem{remark}[theorem]{Remark}

\newcommand{\leftexp}[2]{{\vphantom{#2}}^{{\rm #1}}{#2}}
\newcommand{\leftbase}[2]{{\vphantom{#2}}_{{\rm #1}}{#2}}

\newcommand{\A}{\mathcal{A}}
\newcommand{\D}{\mathcal{D}}
\def\d{\partial}

\def\la{\lambda}
\def\si{\sigma}
\def\Si{\Sigma}
\def\ge{\epsilon}

\def\La{\Lambda}

\newcommand{\q}{\mathbf{q}}

\newcommand{\LD}{\Big\langle}
\newcommand{\RD}{\Big\rangle}

\newcommand{\be}{\mathbf{e}}

\makeatletter

\newcommand{\Rmnum}[1]{\expandafter\@slowromancap\romannumeral #1@}
\makeatother

\newtheorem{thm}{Theorem}[section]

\newtheorem{conj}[thm]{Conjecture}

\theoremstyle{definition}

\newcommand{\beq}{\begin{equation}}
\newcommand{\eeq}{\end{equation}}
\newcommand{\bqa}{\begin{eqnarray}}
\newcommand{\eqa}{\end{eqnarray}}
\newcommand{\beqa}{\begin{equation*}}
\newcommand{\ben}{\begin{eqnarray*}}
\newcommand{\eeqa}{\end{equation*}}
\newcommand{\een}{\end{eqnarray*}}

\newcommand{\B}{\mathcal{B}}
\newcommand{\C}{\mathbb{C}}

\newcommand{\F}{\mathcal{F}}
\renewcommand{\H}{\mathcal{H}}

\renewcommand{\L}{\mathcal{L}}
\newcommand{\M}{\mathcal{M}}

\renewcommand{\O}{\mathcal{O}}
\renewcommand{\P}{\mathbb{P}}

\renewcommand{\S}{\mathcal{S}}

\newcommand{\T}{\mathcal{T}}

\newcommand{\Z}{\mathbb{Z}}

\newcommand{\s}{\mathbf{s}}

\newcommand{\f}{\mathbf{f}}
\renewcommand{\r}{\mathbf{r}}
\renewcommand{\a}{\mathbf{a}}
\newcommand{\m}{\mathbf{m}}
\renewcommand{\v}{\mathbf{v}}

\renewcommand{\t}{\mathbf{t}}
\newcommand{\x}{\mathbf{x}}

\title[Global mirror symmetry for ISES]{Global mirror symmetry for invertible simple elliptic singularities}
\author{Todor Milanov}
\author{Yefeng Shen}
\address{Kavli Institute for the Physics and Mathematics of the Universe (WPI),
Todai Institutes for Advanced Study, The University of Tokyo, Kashiwa, Chiba 277-8583, Japan}
\email{todor.milanov@ipmu.jp}

\address{Kavli Institute for the Physics and Mathematics of the Universe (WPI),
Todai Institutes for Advanced Study, The University of Tokyo, Kashiwa, Chiba 277-8583, Japan}
\email{yefeng.shen@ipmu.jp}

\begin{document}

\maketitle

\tableofcontents
\addtocontents{toc}{\protect\setcounter{tocdepth}{1}}

\section{Introduction}

In the famous mirror symmetry paper \cite{CDGP}, the authors described
a duality of Calabi-Yau 3-folds that exchanges the A-model with the
B-model. The A-model contains information such as K\"ahler structure
and Gromov-Witten invariants while the B-model contains information
such as complex structures and periods integrals. However, this picture is not complete since the complex moduli usually has a nontrivial topology while the K\"ahler moduli does not. Consider the quintic 3-fold as an example. The mirror is a family of quintic 3-folds
$$\sum_{i=1}^{5}X_i^5-5\psi\prod_{i=1}^{5}X_i=0$$
quotient by $(\mathbb{Z}/5\mathbb{Z})^3$, with $\psi\in\mathbb{P}^1$. This new family contain special limits $\psi=0,\infty$ and fifth roots of unity, which are referred to as the Gepner point, the large complex structure limit point and the conifold limits. 
The mirror theorem asserts that the contractible K\"ahler moduli of quintic 3-fold is mirror to a neighborhood of the large complex structure limit \cite{G0}, \cite{LLY}. 

On the other hand, it is implicit in physics that we should study the entire complex moduli and all the special limits. This global point of view leads to BCOV-holomorphic anomaly equation \cite{BCOV} and recent spectacular physics predictions of the Gromov-Witten invariants of quintic 3-fold up to genus 52 \cite{HKQ}. 

Landau-Ginzburg phases are introduced as part of the global picture, to describe the neighborhood of the Gepner point, or its mirror. Recently, a candidate of Landau-Ginzburg A-model has been constructed by Fan, Jarvis and Ruan based on a proposal of Witten \cite{FJR}, \cite{FJR1}. It is now called the Fan-Jarvis-Ruan-Witten theory (FJRW theory). It is a Gromov-Witten type theory which counts solutions of Witten equations. Based on this construction, Ruan proposed a mathematical formulation of Landau-Ginzburg/Calabi-Yau (LG/CY) correspondence \cite{Ru1}. This connects the FJRW theory and Gromov-Witten theory for a pair of same initial data. In \cite{CR1}, Chiodo and Ruan addressed the idea of global mirror symmetry to build a bridge for LG/CY correspondence. In short, in this picture, the FJRW theory is formulated as the mirror theory for the Gepner point.

\subsection{The LG/CY correspondence via global mirror symmetry}\label{gms-int}

Let us briefly recall the general setup for the LG/CY correspondence.
Recall that a polynomial $W$ is called {\em quasi-homogeneous}  if there are rational weights $q_i$ for each $X_i$, such that
$$
W(\lambda^{q_1}X_1, \dots, \lambda^{q_N}X_N)=\lambda W(X_1,\dots, X_N), \quad \forall\lambda\in \C^*.
$$
The polynomial $W$ is called {\em non-degenerate} if: (1) $W$
has isolated critical point at the origin; (2) the choice of all $q_i\in(0,\frac{1}{2}]$
is unique. Let $W(\x)$ be a quasi-homogeneous non-degenerate polynomial,
$$
W(\x)=
\sum_{i=1}^{s} \prod_{j=1}^{N}X_{j}^{a_{ij}},\quad \x=(X_1,\dots,X_N).
$$
We say that $W$ is 
\emph{invertible} if its {\em exponent matrix} $E_W=\left(a_{ij}\right)_{s\times N}$ is an invertible matrix.
A diagonal matrix ${\rm diag}(\lambda_1, \dots, \lambda_N)$ is called a {\em diagonal symmetry} of $W$ if
\ben
W(\la_1 X_1, \dots, \la_N X_N)=W(X_1, \dots, X_N), \quad \la_i\in\C^*.
\een
Let $G_{W}$ be {\em the group of all diagonal symmetries} of $W$. It contains an element
$$
J_W={\rm diag}\left(\exp(2\pi\sqrt{-1} q_1), \dots, \exp(2\pi\sqrt{-1} q_N)\right).
$$
If the {\em Calabi-Yau condition $(\sum_i q_i=1$)} holds,
$X_W=\{W=0\}$ is a Calabi-Yau hypersurface in the weighted
projective space $\mathbb{P}^{N-1}(c_1, \dots, c_N)$, where
$q_i=c_i/d$ for a common denominator $d$. The element $J_W$ acts
trivially on $X_W$, while for any group $G$ such that $\langle
J_W\rangle \subseteq G\subseteq G_{W}$, the group
$\widetilde{G}=G/\langle J_W\rangle$ acts faithfully on $X_W$.  The
LG/CY correspondence \cite{Ru1} predicts that the ancestor potential ($\mathscr{A}^{\rm FJRW}_{W,G}$) of the FJRW theory for $(W, G)$ is the same as the total ancestor potential ($\mathscr{A}^{\rm GW}_{\mathcal{X}}$) of the GW theory for $\mathcal{X}=X_{W}/\widetilde{G}$, up to
analytic continuation and the quantization of a symplectic
transformation. Both $\mathscr{A}^{\rm FJRW}_{W,G}$ and $\mathscr{A}^{\rm GW}_{\mathcal{X}}$ will be defined in Section \ref{sec:3}.

For an invertible polynomial $W$, its \emph{transpose} $W^T$ is the unique invertible polynomial such that $E_{W^T}=(E_W)^T$, where $(E_W)^T$ is the transpose matrix of $E_W$.
The role of the transpose $W^T$ in mirror symmetry was first studied in \cite{BH} by Berglund and
H\"{u}bsch. Later, Krawitz introduced a mirror group $G^T$\cite{K}. Now a pair $(W^T, G^T)$ is referred to as the Berglund-H\"{u}bsch-Krawitz mirror (BHK mirror) of a pair $(W, G)$.

In order to describe the analytic continuation in the LG/CY correspondence for the pair $(W,G)$, Chiodo and Ruan \cite{CR1} addressed the idea of global mirror symmetry. They proposed to consider a global LG B-model for the BHK mirror $(W^T,G^T)$. Such a global moduli contains a Gepner point and a large complex structure limit point. 
Then the FJRW theory is formulated as the mirror theory for the Gepner point, and the GW theory is formulated as the mirror theory for the large complex structure limit point. The LG/CY correspondence is obtained by connecting the Gepner point and the large complex structure limit point on the global moduli.
This works extremely well for $G=G_W$. In this case, the mirror group $G^{T}$ is the trivial group and the Saito-Givental theory of $W^T$ is expected to be the right object of the global LG B-model. 
If $G\neq G_W$, a global Calabi-Yau B-model \cite{CR2}, \cite{CIR} is used to replace the LG B-model for the genus zero theory. However, a mathematical theory for the higher genus of such a global B-model is still not available.  On the other hand, Costello and Li have a different approach to construct a higher genus on the special limits for both CY B-model and LG B-model \cite{CL}, \cite{L}. 

\subsection{Special limits in Saito-Givental theory}
In this paper, we will study the special limits (see Definition
\ref{def-special-limit} below) in the Saito-Givental theory of a one-parameter
family deformation of an invertible simple elliptic singularities (ISES for brevity) and their geometric mirrors. 
All ISESs are listed in Table \ref{table-1}. They are classified into
three different types $E_{\mu-2}^{(1,1)}, \mu=8,9,10$, depending on
the Milnor number $\mu$. 
\begin{table}[h]\center
\caption{Invertible simple elliptic singularities}
\label{table-1}
\begin{tabular}{|l|l|l|l|}
\hline
&$E_6^{(1,1)}$&$E_7^{(1,1)}$&$E_8^{(1,1)}$\\
\hline
Fermat&$X_1^3+X_2^3+X_3^3$&$X_1^4+X_2^4+X_3^2$&$X_1^6+X_2^3+X_3^2$\\
\hline
Fermat+Chain&$X_1^2X_2+X_2^3+X_3^3$&$X_1^3X_2+X_2^4+X_3^2$&$X_1^4X_2+X_2^3+X_3^2$\\
&&$X_1^2X_2+X_2^2+X_3^4$&$X_1^3X_2+X_2^2+X_3^3$\\
\hline
Fermat+Loop&$X_1^2X_2+X_1X_2^2+X_3^3$&$X_1^3X_2+X_1X_2^3+X_3^2$&\\
\hline
Chain&$X_1^2X_2+X_2^2X_3+X_3^3$&$X_1^3X_2+X_2^2X_3+X_3^2$&\\
\hline
Loop&$X_1^2X_2+X_2^2X_3+X_1X_3^2$&&\\
\hline
\end{tabular}
\end{table}

Let $W$ be an invertible simple elliptic singularity (ISES) as in Table \ref{table-1}. Saito's theory of
primitive forms yields a semisimple Frobenius structure on the
miniversal deformation space $\S$ of $W$ \cite{S0}. 
Using Givental's higher-genus reconstruction formalism \cite{G1}, \cite{G2},
we define the {\em total ancestor potential} $\mathscr{A}_{W}^{\rm
  SG}(\s)$ for every semisimple point $\s\in\S$. Details of the Saito-Givental theory will be introduced in
Section \ref{sec:2}.   

Let $\phi_{-1}$ be a weighted-homogeneous monomial of degree 1 that
represents a non-zero element in the Jacobi algebra $\mathscr{Q}_W$ of
$W$. Following the terminology in the physics literature we refer to
$\phi_{-1}$ and to the corresponding deformations $W_\si=W+\si \phi_{-1}$  
respectively as a {\em marginal monomial} and a {\em marginal
  deformation}. The space of marginal deformations $\Si$ consists of
all $\si\in \C$, s.t., the polynomial $W_\si(X)$ has an isolated
critical point at $X=0$. A miniversal deformation of $W$ can be
constructed in such a way that $\S=\Si\times \C^{\mu-1}$.  In the settings of
singularity theory $\mathscr{A}_{W}^{\rm  SG}(\s)$ extends
analytically for all $\s\in \mathcal{S}$ (see \cite{CI,KS,MR} for the
case of simple-elliptic singularity and \cite{M} for the general
case). Hence, there is a unique limit 
$$\mathscr{A}_W^{\rm SG}(\si)=\lim_{\s\to(\si,{\bf
    0})}\mathscr{A}_{W}^{\rm SG}(\s)$$
The space of marginal deformations $\Si$ is a
punctured plane, i.e., $\Si=\C\setminus{\{p_1,\dots,p_l\}}$. An element ${\rm diag}(\la_1,\cdots,\la_n)\in G_W$ acts on $\Si$ by $\si\mapsto \left(\prod_{i=1}^N \la_i\right)\si$.
We say that $\si\in \Si$ is an {\em orbifold point} if it has a
non-trivial stabilizer in $Z:=G_W/({\rm SL}_N(\mathbb{C})\cap G_W)$. Follow Remark 4.2.5 in Chiodo--Ruan \cite{CR1}, $Z$ is a cyclic group and only $\si=0$ is an orbifold point.{}
\begin{definition}\label{def-special-limit}
We call $\si\in\Sigma\cup\{\infty\}$ a {\em special limit point} if 
it is a puncture $\si=p_i$, an orbifold point, or  $\si=\infty$.
Two special limit points are isomorphic if we can identify the corresponding
limits of $\mathscr{A}^{\rm SG}(\sigma)$. We say that
$\si$ is 
\begin{itemize}
\item[a)]
an FJRW-point (or a $(W',G')$-FJRW point) if there exists a pair $(W',G')$, such that $\mathscr{A}_W^{\rm SG}(\si)=\mathscr{A}^{\rm FJRW}_{W',G'},$ where $\mathscr{A}^{\rm FJRW}_{W',G'}$ is the total ancestor potential of the FJRW theory for $(W',G')$.
\item[b)]
a GW-point (or a $\mathcal{X}$-GW point) if there exists an orbifold $\mathcal{X}$,
such that $\mathscr{A}_W^{\rm SG}(\si)=\mathscr{A}^{\rm GW}_{\mathcal{X}}$, where $\mathscr{A}^{\rm GW}_{\mathcal{X}}$ is the total ancestor potential of $\mathcal{X}$.
\end{itemize}
\end{definition}
In the physics literature, the FJRW-points are known as {\em
  Gepner points}, while the GW-points are also known as {\em large
  complex structure limit} points. 
Let $E_{\si}$ be the elliptic curve in $\mathbb{P}^2(c_1,c_2,c_3)$,
defined by $W_{\si}=0$. Let $j(\sigma)$ be the $j$-invariant of
$E_{\sigma}$ and $\mu$ be the Milnor number of $W$. Based on the
calculations in \cite{MR,KS}, we propose the following conjecture to
understand the mirror symmetry and to classify the special limit
points for invertible simple elliptic singularities.  
\begin{conj}\label{main-conj}
In the set of all ISESs with a given marginal deformation, we have
\begin{itemize}
\item[a)] Each of the special limit points is either an FJRW point
  corresponding to some ISES or a GW point corresponding to some
  elliptic orbifold $\mathbb{P}^1$. 

\item[b)] The map that assigns to each point $\si$ the pair $(\mu,j(\si))$
  induces a one-to-one correspondence between the isomorphism classes of FJRW and GW
  points and the set of pairs $(\mu,j)$, with $\mu=8,9,10$ and $j=0,1728,\infty$.
\end{itemize}
\end{conj}

\begin{remark}
In the quintic case, the points $p_1,\dots,p_l$
are neither FJRW or GW points. Usually they are referred to as conifold points. It is still not known
if an appropriate geometric mirror model exists for a conifold point. If
Conjecture \ref{main-conj} holds, then there are no conifold points for
ISESs.
\end{remark}

Our first main result is that Conjecture \ref{main-conj} is true for $\si=0$.
\begin{theorem}\label{t1} 
Let $W$ be an invertible polynomial of type
$E_{\mu-2}^{(1,1)}$.
\begin{itemize}
\item[a)]
If $W^T$ belongs to Tables \ref{ring-E6}, \ref{ring-E7}, or \ref{ring-E8}, then $\sigma=0$ is a $(W^T,G_{W^T})$-FJRW point.
\item[b)] 
$\si=0$ is always an FJRW point. The isomorphism class of $\si=0$ as a FJRW point is determined by
$(\mu,j(\si))$, with $\mu=8,9,10$, $j(\si)=0$ or $1728$.  
\end{itemize}
\end{theorem}
Note that in Theorem \ref{t1} we have excluded some of the
polynomials appearing in Table \ref{table-1}. This is because we do not know how to
compute all the FJRW invariants for them(see Section \ref{sec-3.2}). 
\begin{remark}
Part a) of Theorem \ref{t1} is a particular case of the so called {\em LG-LG
  mirror symmetry} \cite{CR1}. 
In the terminology of Chiodo-Ruan, the
pair  $(W,\{1\})$ is called the BHK-mirror of $(W^T,G_{W^T})$.
\end{remark}

Our next main result is the following.
\begin{theorem}\label{t2}
Conjecture \ref{main-conj} is true for the ISESs of Fermat type.
\end{theorem}
A complete answer to Conjecture \ref{main-conj} for all special limits
$\si\neq0$ of all ISESs can be achieved in a similar fashion. However,
our approach works on a case by case basis. Verifying the conjecture
in the remaining cases is a work in progress \cite{LMS}. It will be
nice if one can find a more conceptual prove of Conjecture
\ref{main-conj} that does not rely on the explicit form of the
polynomial $W$. 
\begin{remark}
A corollary of Theorem \ref{t2} is that in the case of Fermat
$E^{(1,1)}_8$, the special limit point $\si=\infty$ is a FJRW-point and
all other special limit points $\si\neq0$ are GW-points.  
This is another surprising result. Usually  in a one-parameter B-model
family \cite{DM}, the point $\sigma=\infty$ is expected to be a
GW-point.
\end{remark}

As a corollary, we get various correspondences of LG/LG-type or of LG/CY-type. 
\begin{corollary}
The total ancestor potentials of the GW theories and the FJRW theories
that are obtained as special limits of the total ancestor potential of
a given Fermat type ISES are related by analytic continuation and 
quantizations of symplectic transformations.
\end{corollary}

In \cite{HKQ}, the authors use a certain gap condition at the conifold points
and regularity at Gepner points to compute the GW invariants of the
quintic 3-fold up to genus 52. In our case, it is much nicer since
Theorem \ref{t2} implies the following corollary. 
\begin{corollary}
Let $W$ be a Fermat simple elliptic singularity and $g\geq 0$ be an
arbitrary integer. The
genus-$g$ Saito-Givental correlation functions are holomorphic near
the special limit points. 
\end{corollary}
Using the global B-model and mirror symmetry, Milanov--Ruan \cite{MR}
proved that the GW invariants of any genus of the orbifolds
$\mathbb{P}^1_{3,3,3}$,$\mathbb{P}^1_{4,4,2}$ and
$\mathbb{P}^1_{6,3,2}$ are quasi-modular forms for some modular group
$\Gamma(W)$.  
We will compute $\Gamma(W)$ for the
Fermat polynomials $W$ of type $E_{\mu-2}^{(1,1)},$ $\mu=9,10$ in a subsequent paper \cite{MS}. 
It would be interesting to see whether this helps to express the higher genus GW invariants explicitly in  closed forms, as polynomials of ring generators of quasi-modular forms.

\subsection{Plan of the paper}

The paper is organized as follows.
In Section \ref{sec:2}, we recall 
Givental's construction of the B-model Gromov-Witten type
potential in the setting of Saito's theory of primitive forms. We
also derive the Picard-Fuchs equations satisfied by the various period
integrals. In Section \ref{sec:3}, we
discuss the two types of geometric theories:
the Gromov-Witten theory of elliptic orbifold lines and the FJRW
theory of simple elliptic singularities. We also recall the
reconstruction theorem in both theories. In
Sections \ref{MS-0}, we establish the LG-LG mirror symmetry for $\si=0$ (Theorem \ref{t1}) by comparing the B-model
constructed in Section \ref{sec:2} and the FJRW A-models constructed in
Section \ref{sec:3}.  In Section \ref{MS-punctures}, 
we establish the global mirror symmetry  for Fermat polynomials by proving Theorem \ref{t2}.

\subsection{Acknowledgement}
We thank Yongbin Ruan for his insight and support for this project.
Both authors would like to thank Kentaro Hori, Hiroshi Iritani, and
Kyoji Saito for many stimulating conversations. The first author benefited from
conversations with Satoshi Kondo and Charles Siegel, while the second
author would like to thank Alessandro Chiodo, Igor Dolgachev, Huijun Fan, Tyler
Jarvis, Jeffrey Lagarias, Sijun Liu, Noriko Yui and Jie Zhou for helpful discussions. We thank Arthur
Greenspoon for editorial assistance. The first authors is
supported by Grant-In-Aid and by the World Premier International
Research Center Initiative (WPI Initiative), MEXT, Japan.

\section{Global B-model for simple elliptic singularities}\label{sec:2}

Let $W$ be an invertible polynomial from Table \ref{table-1}. We would like to recall Saito's
theory of primitive forms which yields a Frobenius structure on
$\S$. Following Givental's higher genus reconstruction formalism we
will introduce the total ancestor potential of $W$. We also
derive a system of hypergeometric equations that determines the
restriction of the flat coordinates of the Frobenius manifold $\S$ to
$\Si$. 

\subsection{Miniversal deformations}
In this paper, we are interested in invertible polynomials with 3 variables and the weights
$q_i=1/a_i(1\leq i\leq 3)$ for some positive integers $a_i$ satisfying the Calabi-Yau
condition
$$
(a_1,a_2,a_3)=(3,3,3),\ (4,4,2),\ \mbox{and}\ (6,3,2).
$$
We denote the corresponding classes of invertible polynomials
respectively by $E_6^{(1,1)}$, $E_7^{(1,1)}$, and
$E_8^{(1,1)}$. Modulo permutation of the variables $X_i(1\leq i\leq
N)$ there are 13 types of invertible polynomials (see Table
\ref{table-1}). We refer to these polynomials as invertible simple
elliptic singularities.  Let $\mathscr{Q}_{W}$
be the \emph{Jacobian algebra} of $W$,
$$
\mathscr{Q}_{W}=\C[X_1,X_2,X_3]/(\d_{X_1}W,\d_{X_2}W,\d_{X_3}W).
$$
Let us fix a set $\mathfrak{R}$ of
weighted homogeneous monomials 
\beq\label{JA:basis}
\phi_\r(\x)=X_1^{r_1}X_2^{r_2}X_3^{r_3},\quad \r=(r_1,r_2,r_3)\in\mathfrak{R},
\eeq
such that their projections in $\mathscr{Q}_W$ form a basis. The dimension
of $\mathscr{Q}_{W}$
called the {\em multiplicity} of the critical point or {\em Milnor number}
and it will be denoted by $\mu$. There
is precisely one monomial of top degree, say
$\phi_{\m}$, $\m=(m_1,m_2,m_3)\in\mathfrak{R}$. We fix a deformation of $W$
of the following form: 
\beq\label{ises}
W_\si(\x)=W(\x)+\si\, \phi_\m (\x),\quad \si\in \Si,
\eeq
where $\Sigma\subset \C$ is the set
of all $\si\in \C$ such that $W_{\si}(\x)$ has only isolated critical
points. Such deformations do not change the multiplicity of the
critical point at $\x=0$. The polynomials \eqref{ises} are families of simple elliptic singularities of type
$E_{\mu-2}^{(1,1)}$ (see \cite{S2}).
More generally, we consider a miniversal deformation (see
e.g. \cite{AGV}) of $W$
\beq\label{min-def}
F(\s,\x)=W(\x)+\sum_{\r\in \mathfrak{R}} s_\r\, \phi_\r (\x).
\eeq
It is convenient to adopt two notations for the deformation
parameters. Namely, put
\ben
\s=\{s_\r\}_{\r\in \mathfrak{R}} = (s_{-1},s_0,s_1,\dots,s_{\mu-2}),
\een
where the second equality is obtained by putting an order on the
elements $\r\in \mathfrak{R}$ and enumerating them with the integers
from $-1$ to $\mu-2$ in such a way that 
\ben
s_{-1}=s_\m=\si,\quad s_0=s_{{\bf 0}},\quad {\bf 0}=(0,0,0)\in \mathfrak{R}.
\een
The moduli space of miniversal deformations, i.e., the range of the
parameters $s_\r$ is then defined to be the affine space $\S=\Si\times
\C^{\mu-1}$. 
Furthermore, each $s_\r$ is assigned a degree so that $F(\s,\x)$ is weighted-homogeneous of
degree 1. Note that the parameter $s_{\m}=\si$ has degree $0$. 
Following the terminology in physics we call $s_{\m}$ and
 $\phi_{\m}$ {\em marginal}. Note that $W_\si(\x)$ is the restriction
 of $F(\s,\x)$ to the subspace $\Si$ of marginal deformations. 
 Except for Fermat case, there is more than one choice
 of a marginal monomial. For example, $X_1X_2X_3$, $X_1^4X_3$
 are both marginal for $W=X_1^3X_2+X_2^2+X_3^3$.

\subsection{Saito's theory}
Let $C$ be the critical variety of the miniversal deformation $F(\s,\x)$ (see
(\ref{min-def})), i.e., the support of the sheaf
\ben
\O_C:=\O_X/\langle \d_{X_1} F,\d_{X_2}F,\d_{X_3}F\rangle,
\een
where $X=\S\times \C^3$. Let $q:X\to \S$ be the projection on the
first factor. The Kodaira--Spencer map ($\T_\S$ is the sheaf of
holomorphic vector fields on $\S$) 
\ben
\T_\S\longrightarrow q_*\O_C,\quad \d/\d s_i\mapsto \d F/\d s_i\ {\rm
  mod} \ (F_{X_1},F_{X_2},F_{X_3})
\een 
is an isomorphism, which implies that for any ${\bf s}\in\S,$ the tangent space $T_{\bf s}\S$
is equipped with an associative commutative multiplication
$\bullet_{\bf s}$ depending holomorphically on $\bf{s}\in \S$. If in
addition we have a volume form $\omega=g({\bf s},{\bf x})d^3{\bf x},$
where $d^3{\bf x}=dX_1\wedge dX_2\wedge dX_3$ is the standard volume form, then $q_*\O_C$ (hence $\T_\S$ as well) is equipped with the {\em residue pairing}:
\beq\label{res:pairing}
\LD\psi_1,\psi_2\RD= \frac{1}{(2\pi i)^3} \ \int_{\Gamma_\ge} \frac{\psi_1({\bf s,y})\psi_2({\bf s,y})}{F_{y_1}F_{y_2} F_{y_3}}\, \omega,
\eeq
where ${\bf y}=(y_1,y_2,y_3)$ is a unimodular coordinate system for the volume form, i.e., $\omega=d^3{\bf y}$, and $\Gamma_\ge$ is a real $3$-dimensional cycle supported on $|F_{X_i}|=\ge$ for $1\leq i\leq 3.$

Given a semi-infinite cycle
\beq\label{cycle}
\A\in  \lim_{ \longleftarrow } H_3(\C^3, (\C^3)_{-m} ;\C)\cong \C^\mu,
\eeq
where
\beq\label{level:lower}
(\C^3)_m=\{{\bf x}\in \C^3\ |\ {\rm Re}(F({\bf s,x})/z)\leq m\}.
\eeq
Let $d_\S$ be the de Rham differential on $\S$. Put
\beq\label{osc_integral}
J_{\A}({\bf s},z) = (-2\pi z)^{-3/2} \, zd_\S \, \int_{\A} e^{F({\bf s,x})/z}\omega,
\eeq
The oscillatory integrals $J_\A$ are by definition sections of the cotangent sheaf $\T_\S^*$.

According to Saito's theory of primitive forms \cite{S0}, there exists a volume form $\omega$ such that the residue pairing is flat and the oscillatory integrals satisfy a system of differential equations, which have the form
\beq\label{frob_eq1}
z\d_i J_\A({\bf t},z) = \d_i \bullet_{\bf t} J_\A({\bf t},z),\quad \d_i:=\d/\d t_i\ (-1\leq i\leq\mu-2),
\eeq
in flat-homogeneous coordinates ${\bf t}=(t_{-1},t_0,\dots,t_{\mu-2})$ and the multiplication is
defined by identifying vectors and covectors via the residue pairing. 
Using the residue pairing, the flat structure, and the
Kodaira--Spencer isomorphism we have:
\ben
T^*\S\cong T\S\cong \S\times T_0\S\cong \S\times \mathscr{Q}_{W}.
\een
Due to homogeneity, the integrals satisfy a differential equation:
\beq
\label{frob_eq2}
(z\d_z + E)J_{\A}({\bf t},z) =  \Theta\, J_\A({\bf t},z), \quad z\in \C^*
\eeq
where $E$ is the {\em Euler vector field} 
\ben
E=\sum_{i=-1}^{\mu-2}d_it_i \d_i,\quad (d_i:={\rm deg}\, t_i={\rm deg}\, s_i),
\een
and $\Theta$ is the so-called {\em Hodge grading operator }
\ben
\Theta:\T^*_S\rightarrow \T^*_S,\quad \Theta(dt_i)=\left(\frac{1}{2}-d_i\right)dt_i.
\een
The compatibility of the system \eqref{frob_eq1}--\eqref{frob_eq2}
implies that the residue pairing, the multiplication, and the Euler
vector field give rise to a {\em conformal Frobenius structure} of
conformal  dimension $1$. We refer to B. Dubrovin \cite{Du} for the
definition and more details on Frobenius structures and to C. Hertling
\cite{He} or to Atsushi--Saito 
\cite{SaT} for more details on constructing a Frobenius structure
from a primitive form.

\subsection{The primitive forms}
The classification of primitive forms in general is a very difficult
problem. In the case of simple elliptic singularities however, all
primitive forms are known (see \cite{S0}). They are given by  $\omega=d^3\textbf{x}/\pi_{A}(\sigma)$,
where $\pi_A(\si)$ is the period (\ref{pi_A}). As we will prove below,
these periods are solutions to the hypergeometric equation \eqref{PF:1},
so a primitive form may be equivalently fixed by fixing a solution to
the differential equation that does not vanish on $\Si$. Note that
since $\pi_A(\si)$ is multi-valued function, the corresponding
Frobenius structure on $\S$ is multi-valued as well. In other words,
the primitive form gives rise to a Frobenius structure on the
universal cover $\widetilde{\S}\cong \mathbb{H}\times \C^{\mu-1}$. 

The key to the primitive form is the Picard-Fuchs differential
equation for the periods of the so-called {\em elliptic curve at infinity}
\beq\label{ef}
E_{\sigma}:=\Big\{[X_1:X_2:X_3]\in\mathbb{CP}^2(c_1,c_2,c_3)\,\Big\vert\,W_\si=0\Big\},
\eeq
where $c_i=d/a_i$, $1\leq i\leq 3$ and $d$ is the least common
multiple of $a_1,a_2$, and $a_3$. Note that $E_\si$ are the fibers of
an elliptic fibration over $\mathbb{CP}^1=\C\cup\{\infty\}$ whose non-singular fibers
are parametrized by $\Si\subset \C\subset \mathbb{CP}^1.$
Note that $\mathrm{Res}_{E_{\sigma}}\Omega$, where
$$
\Omega:=\frac{dX_1\wedge dX_2\wedge dX_3}{dW_{\sigma}},
$$
is a Calabi-Yau form of the elliptic curve $E_{\sigma}$. For every
$A\in H_{1}(E_{\sigma}),$ we define 
\begin{equation}\label{pi_A}
\pi_{A}(\sigma)=\int_{A}\mathrm{Res}_{E_{\sigma}}\Omega.
\end{equation}
It is well known that the period integrals are solutions to a Fuchsian
differential equation. In particular, we obtain the following lemma, 
\begin{lemma}\label{deq}
Let $\delta=\si{\d}/{\d \si}$, the elliptic period \eqref{pi_A} satisfies the Picard-Fuchs equation 
\beq\label{PF:1}
\delta(\delta-1)\ \pi_A(\si)=C\, \sigma^l(\delta+l\alpha)(\delta+l\beta)\ \pi_A(\si), \quad \alpha+\beta=1-\frac{1}{l}, \quad C=\prod_{i=1}^3 \left(-\frac{l_i}{l}\right)^{l_i}.
\eeq
If we put $x=C\, \sigma^l$, $\gamma=\alpha+\beta$, the equation \eqref{pi_A} becomes a hypergeometric equation
\beq\label{HG:1}
x(1-x)\frac{d^2\pi_A}{dx^2}+\Big(\gamma-(1+\alpha+\beta)x\Big)\frac{d\pi_A}{dx}-(\alpha\beta)\,\pi_A=0
\eeq
We call $(\alpha,\beta,\gamma)$ the weights system of the hypergeometric equation. 
The explicit values are listed in Tables
\ref{PF-E6:weights}, \ref{PF-E7:weights} and \ref{PF-E8:weights} below.
\end{lemma}
In particular,  $\Sigma=\C\backslash \{p_1,\dots,p_l\}$, where $p_i$ are 
the singularities of the Picard--Fuchs equation \eqref{PF:1}. All the singular points are
\beq\label{sing-pts}
p_i=C^{-1/l}\eta^i, \quad 1\leq i\leq l,\quad \eta=\exp(2\pi\sqrt{-1}/l),
\eeq
For our purposes, We will give a proof of this lemma in Section \ref{PF:system} following the approach of S. G\"ahrs (see \cite{G}). To find out $\alpha,\beta$ and $\gamma$, we will need {\em the mirror weight} $q_i^T$, which is the weight of $X_i$ in the BHK mirror $W^T$ and {\em a charge vector} 
$(l_1,l_2,l_3,-l)\in \Z^4$
by choosing the minimal
$l\in\mathbb{Z}_{>0}$ such that
\beq\label{weights}
(l_1,l_2,l_3)= l\,\m\, E_{W}^{-1},\quad \m=(m_1,m_2,m_3).
\eeq

\begin{table}[h]\center
\caption{$E_{6}^{(1,1)}$}\label{PF-E6:weights}
\begin{tabular}{|c|c|c|c|c|}
\hline
$W$ & $m_1,m_2,m_3$&$l_1,l_2,l_3,l$& $q_1^T,q_2^T,q_3^T$&$\alpha,\beta,\gamma$ \\
\hline
  $X_1^3+X_2^3+X_3^3$
  &$1,1,1$
  &$1,1,1,3$
    &$\frac{1}{3},\frac{1}{3},\frac{1}{3}$
  &$\frac{1}{3},\frac{1}{3},\frac{2}{3}$
  \\

  \hline
  $X_1^2X_2+X_2^3+X_3^3$
  &$2,0,1$
  &$3,-1,1,3$
    &$\frac{1}{2},\frac{1}{6},\frac{1}{3}$
  &$\frac{1}{6},\frac{1}{2},\frac{2}{3}$
  \\

  \hline
  $X_1^2X_2+X_2^3+X_3^3$
  &$0,2,1$
  &$0,2,1,3$
    &$\frac{1}{2},\frac{1}{6},\frac{1}{3}$
  &$\frac{1}{12},\frac{7}{12},\frac{2}{3}$
  \\

  \hline
  $X_1^2X_2+X_1X_2^2+X_3^3$
  &$1,1,1$
  &$1,1,1,3$
    &$\frac{1}{3},\frac{1}{3},\frac{1}{3}$
  &$\frac{1}{3},\frac{1}{3},\frac{2}{3}$
  \\

  \hline
  $X_1^2X_2+X_1X_2^2+X_3^3$
  &$2,0,1$
  &$4,-2,1,,3$
    &$\frac{1}{3},\frac{1}{3},\frac{1}{3}$
  &$\frac{1}{12},\frac{7}{12},\frac{2}{3}$
  \\

  \hline
  $X_1^2X_2+X_2^2X_3+X_3^3$
  &$2,0,1$
  &$2,-1,1,2$
    &$\frac{1}{2},\frac{1}{4},\frac{1}{4}$
  &$\frac{1}{4},\frac{1}{4},\frac{1}{2}$
  \\

  \hline
  $X_1^2X_2+X_2^2X_3+X_3^3$
  &$0,3,0$
  &$0,3,-1,2$
    &$\frac{1}{2},\frac{1}{4},\frac{1}{4}$
  &$\frac{1}{12},\frac{5}{12},\frac{1}{2}$
  \\

  \hline
  $X_1^2X_2+X_2^2X_3+X_3^3$
  &$0,1,2$
  &$0,1,1,2$
    &$\frac{1}{2},\frac{1}{4},\frac{1}{4}$
  &$\frac{1}{4},\frac{1}{4},\frac{1}{2}$
  \\

  \hline
  $X_1^2X_2+X_2^2X_3+X_1X_3^2$
  &$1,1,1$
  &$1,1,1,3$
    &$\frac{1}{3},\frac{1}{3},\frac{1}{3}$
  &$\frac{1}{3},\frac{1}{3},\frac{2}{3}$
  \\

  \hline
  $X_1^2X_2+X_2^2X_3+X_1X_3^2$
  &$3,0,0$
  &$4,-2,1,3$
    &$\frac{1}{3},\frac{1}{3},\frac{1}{3}$
  &$\frac{1}{12},\frac{7}{12},\frac{2}{3}$
  \\
\hline

\end{tabular}
\end{table}

\begin{table}[h]\center
\caption{$E_{7}^{(1,1)}$}\label{PF-E7:weights}
\begin{tabular}{|c|c|c|c|c|}
\hline
$W$ & $m_1,m_2,m_3$&$l_1,l_2,l_3,l$& $q_1^T,q_2^T,q_3^T$&$\alpha,\beta,\gamma$ \\
  \hline
  $X_1^4+X_2^4+X_3^2$
  &$2,2,0$
    &$1,1,0,2$
    &$\frac{1}{4},\frac{1}{4},\frac{1}{2}$
  &$\frac{1}{4},\frac{1}{4},\frac{1}{2}$
  \\

\hline
  $X_1^3X_2+X_2^4+X_3^2$
  &$4,0,0$
  &$4,-1,0,3$
    &$\frac{1}{3},\frac{1}{2},\frac{1}{6}$
  &$\frac{1}{12},\frac{7}{12},\frac{2}{3}$
  \\

\hline
  $X_1^3X_2+X_2^4+X_3^2$
  &$1,3,0$
  &$1,2,0,3$
    &$\frac{1}{3},\frac{1}{2},\frac{1}{6}$
  &$\frac{1}{12},\frac{7}{12},\frac{2}{3}$
  \\

\hline
  $X_1^2X_2+X_2^2+X_3^4$
  &$2,0,2$
  &$2,-1,1,2$
    &$\frac{1}{2},\frac{1}{4},\frac{1}{4}$
  &$\frac{1}{4},\frac{1}{4},\frac{1}{2}$
  \\

\hline
  $X_1^2X_2+X_2^2+X_3^4$
  &$0,1,2$
  &$0,1,1,2$
    &$\frac{1}{2},\frac{1}{4},\frac{1}{4}$
  &$\frac{1}{4},\frac{1}{4},\frac{1}{2}$
  \\

\hline
  $X_1^3X_2+X_1X_2^3+X_3^2$
  &$4,0,0$
  &$3,-1,0,2$
    &$\frac{1}{4},\frac{1}{4},\frac{1}{2}$
  &$\frac{1}{12},\frac{5}{12},\frac{1}{2}$
  \\

\hline
  $X_1^3X_2+X_1X_2^3+X_3^2$
  &$2,2,0$
  &$1,1,0,2$
    &$\frac{1}{4},\frac{1}{4},\frac{1}{2}$
  &$\frac{1}{4},\frac{1}{4},\frac{1}{2}$
  \\

\hline
  $X_1^3X_2+X_2^2X_3+X_3^2$
  &$1,1,1$
  &$1,1,1,3$
    &$\frac{1}{3},\frac{1}{3},\frac{1}{3}$
  &$\frac{1}{3},\frac{1}{3},\frac{2}{3}$
  \\

\hline
  $X_1^3X_2+X_2^2X_3+X_3^2$
  &$1,3,0$
  &$1,4,-2,3$
    &$\frac{1}{3},\frac{1}{3},\frac{1}{3}$
  &$\frac{1}{12},\frac{7}{12},\frac{2}{3}$
  \\

\hline
  $X_1^3X_2+X_2^2X_3+X_3^2$
  &$4,0,0$
  &$4,-2,1,3$
    &$\frac{1}{3},\frac{1}{3},\frac{1}{3}$
  &$\frac{1}{12},\frac{7}{12},\frac{2}{3}$
  \\
\hline
\end{tabular}
\end{table}

\begin{table}[h]\center
\caption{$E_{8}^{(1,1)}$}\label{PF-E8:weights}
\begin{tabular}{|c|c|c|c|c|}
\hline
W & $m_1,m_2,m_3$&$l_1,l_2,l_3,l$& $q_1^T,q_2^T,q_3^T$&$\alpha,\beta,\gamma$ \\
  \hline
  $X_1^6+X_2^3+X_3^2$
  &$4,0,1$
    &$1,2,0,3$
    &$\frac{1}{6},\frac{1}{3},\frac{1}{2}$
  &$\frac{1}{12},\frac{7}{12},\frac{2}{3}$
  \\

\hline
  $X_1^3X_2+X_2^2+X_3^3$
  &$1,1,1$
  &$1,1,1,3$
    &$\frac{1}{3},\frac{1}{3},\frac{1}{3}$
  &$\frac{1}{3},\frac{1}{3},\frac{2}{3}$
  \\

  \hline
  $X_1^3X_2+X_2^2+X_3^3$
  &$4,0,1$
  &$4,-2,1,3$
  &$\frac{1}{3},\frac{1}{3},\frac{1}{3}$
  &$\frac{1}{12},\frac{7}{12},\frac{2}{3}$
  \\

  \hline
  $X_1^4X_2+X_2^3+X_3^2$
  &$2,2,0$
    &$1,1,0,2$
    &$\frac{1}{4},\frac{1}{4},\frac{1}{2}$
    &$\frac{1}{4},\frac{1}{4},\frac{1}{2}$
  \\

  \hline
  $X_1^4X_2+X_2^3+X_3^2$
  &$6,0,0$
    &$3,-1,0,2$
    &$\frac{1}{4},\frac{1}{4},\frac{1}{2}$
  &$\frac{1}{12},\frac{5}{12},\frac{1}{2}$
  \\
\hline
\end{tabular}
\end{table}

\subsection{Picard-Fuchs equations}\label{PF:system}
Let us denote by 
\ben
X_{\s}=\{\x\in \C^3\ |\ F(\s,\x)=1\},\quad \s\in \S.
\een
The points $\s$ for which $X_\s$ is singular form an analytic
hypersurface in $\S$ called the {\em discriminant}. Its complement in $\S$
will be denoted by $\S'$.  We will be interested in the period integrals
\ben
\Phi_{\r}(\s) = \int \phi_{\r}(x)\, \frac{d^3\x}{d F},
\quad
\phi_{\r}(\x)=X_1^{r_1}X_2^{r_2}X_3^{r_3}, 
\quad 
\r=(r_1,r_2,r_3).
\een
They are sections of the middle (or vanishing) cohomology bundle on
$\S'$ formed by $H^2(X_\s,\C)$. Slightly abusing the notation, we denote the restriction to
$s_{-1}=\si$, $s_i=0(0\leq i\leq \mu-2)$ by $\Phi_{\r}(\si)$. 
Following the idea of \cite{G}, we first obtain a GKZ
(Gelfand--Kapranov--Zelevinsky) system of differential equations for
the periods. Using that the period integrals are not polynomial in
$\si$ (they have singularities at the punctures of $\Si$) we can
reduce the GKZ system to a Picard-Fuchs equation.

\subsubsection{The GKZ system}
In order to derive the GKZ system, we slightly modify the
polynomial $W$. By definition $W(\x)=\sum_{i=1}^{3} \phi_{\a_i}(\x),$ where 
$\a_i$ are the rows of the matrix $E_W$. Put
\ben
W_{\v,\si}(\x)=\sum_{i=1}^3 v_i\, \phi_{\a_i}(\x) + \si\,
\phi_{-1}(\x),
\een
where $\v=(v_1,v_2,v_3)$ are some complex
parameters. For simplicity. we omit $\v$ in the notation if $\v=
(1,1,1)$. Let us write $X_\si^{\v,\la}=\{\x\in \C^3\ |\
W_{\v,\si}(\x)=\la\}$. Then we define the period integrals 
\beq\label{GKZ-period}
\Phi_{\r}^{\v,\la}(\si)=\int \phi_{\r}(\x)\, \frac{d^3\x}{dW_{\v,\si}};
\eeq
again one should think that the above integral is a section of the
vanishing cohomology for $W_{\v,\si}(\x).$
The vanishing cohomology bundle is equipped with a Gauss--Manin
connection $\nabla$. The following formulas are well known (see
e.g. \cite{AGV})
\beq\label{GM-con}
\begin{aligned}
\nabla_{\d/\d\la}\ \int \theta  & = \ \ \  \int\frac{d\theta}{d W_{\v,\si}}
\\
\nabla_{\d/\d v_i}\ \int \theta  & = - \int \frac{\d W_{\v,\si}}{\d
  v_i}\, \frac{d\theta}{d W_{\v,\si}} + \int {\rm Lie}_{\d/\d v_i} \theta,
\end{aligned}
\eeq
where $\theta$ is a 2-form on $\C^3$ possibly depending on the
parameters $\v$.
Finally, note that rescaling $X_i\mapsto \la^{q_i} \, X_i (1\leq i\leq 3)$ yields 
\ben
\Phi_{\r}^{\v,\la}(\si) = \lambda^{{\deg\phi}_{\r}}\ \Phi_{\r}^{\v,1}(\si).
\een

Let $\delta_{i}=v_i\d/\d v_i(1\leq i\leq 3)$ and $\delta=\si\d/\d \si$.
\begin{lemma}\label{GKZ-system}
The period integral $\Phi_{\r}^{\v,\la}$ satisfies the system of
differential equations:
\ben
\d_\si^l\prod_{i:l_i<0}\d_{v_i}^{-l_i}\ \Phi & = &
\prod_{i:l_i>0}\d_{v_i}^{l_i}\ \Phi;\\
(\delta_1,\delta_2,\delta_3)\, E_{W}\, \Phi 
+(m_1,m_2,m_3)\delta \,\Phi & = &
-(1+r_1,1+r_2,1+r_3)\Phi.
\een
where the range for $i$ and $j$ in the
first equation on the LHS and the RHS is $0\leq i\leq 3$.
\end{lemma}
\proof
Using \eqref{GM-con} we get the following differential equations:
\ben
\d_{v_i} \ \Phi_{\r}^{\v,\la}\  =-\d_\la \, \Phi_{\r+\a_i}^{\v,\la},\quad
1\leq i\leq 3, 
\een
and
\ben
\d_{\si} \ \Phi_{\r}^{\v,\la}\  =-\d_\la \, \Phi_{\r+\m}^{\v,\la},\quad
\m=(m_1,m_2,m_3),
\een
where $\phi_{\m}(\x)$ is the marginal monomial. 
The first differential equation is equivalent to the identity
$$
l\,m_k-\sum_{i,l_i<0}a_{ik}\,l_i=\sum_{j,l_j>0}a_{jk}\,l_j\,.
$$
which is true by definition (see \eqref{weights}).
For the second equation, using the above formulas we get that the
$i$-th entry on the LHS is 
\ben
-\d_\la\ \int X_i \,\phi_{\r}(X)\, \frac{d^3\x}{dX_i} = -(r_i+1) \int\ \phi_{\r}(X)\, \frac{d^3\x}{dX_i},
\een
where we used formulas \eqref{GM-con} again.
\qed

Let us define the row-vector 
\beq\label{row-vector}
\zeta=(\zeta_1,\zeta_2,\zeta_3)=\r\,
E_W^{-1}.
\eeq 
Note also that the weights $(q_1^T,q_2^T,q_3^T)$ of the
mirror polynomial $W^T$ are precisely 
\beq\label{dual-weights}
(q_1^T,q_2^T,q_3^T)=(1,1,1)\, E_W^{-1}
\eeq
\begin{lemma}\label{GKZ-si}
Let $C=\prod_{i=1}^{3}(-l_i/l)^{l_i}, $ the period integral $\Phi_\r(\si)$ is in the kernel of the differential operator:
\beq\label{semi-GKZ}
\sigma^{-l}\prod_{k=0}^{l-1}\left(\delta-k\right)\prod_{i,l_i<0}\prod_{k=0}^{-l_i-1}\left(\delta+\frac{l\left(q_i^T+\zeta_i+k\right)}{l_i}\right)
-
C\,
\prod_{i,l_i>0}\prod_{k=0}^{l_i-1}\left(\delta+\frac{l\left(q_i^T+\zeta_i+k\right)}{l_i}\right),
\eeq
\end{lemma}
\proof
Using the second equation in Lemma \ref{GKZ-system} we can express the
derivatives $\d_{v_i}=v_i^{-1} \delta_i$ in terms of
$\delta$. Substituting in the first equation we get a higher order
differential equation in $\si$ only. It remains only to notice that the resulting equation is
independent of $v$ and $\la$.
\qed

\subsubsection{Picard-Fuchs equation}
Let  $q_0^T=0, l_0=-l$, and set 
\beq\label{period-ind}
\beta_{i,k}= \frac{1}{l_i}(q_i^T+\zeta_i+k),\quad 0\leq k\leq |l_i|-1.
\eeq
The differential operator in \eqref{semi-GKZ} is the product of a Bessel differential operator
\beq\label{Bessel}
\prod_{i,k} (\delta + l\,\beta_{i,k})
\eeq
and an operator of the form
\beq\label{PF-operator}
\prod_{i',k'} (\delta + l \,\beta_{i',k'}) - C\si^l\prod_{i'',k''}(\delta+l\,\beta_{i'',k''}).
\eeq
This is done by factoring out the common left divisors in the two
summands. There is no pairs $(i',k')$ and $(i'',k'')$ in the
operator \eqref{PF-operator}, such that, $\beta_{i',k'}+1 = \beta_{i'',k''}$. 

\begin{lemma}\label{identity}
The numbers \eqref{period-ind} satisfy the following identity:
$$
\sum_{i:l_i>0}\sum_{k=0}^{l_i-1}\beta_{i,k}-\sum_{0\leq j\leq3:l_j<0}\sum_{k'=0}^{-l_j-1}\left(1+\beta_{j,k'}\right)=\deg\phi_{\r}.$$
\end{lemma}
\proof
By definition 
\begin{eqnarray*}
LHS&=&\sum_{i:l_i>0}\left(\frac{l_i-1}{2}+q_i^T+\zeta_i\right)-\sum_{j;l_j<0}\left(-l_j-\frac{-l_j-1}{2}-q_j^T-\zeta_j\right)-\frac{l-1}{2}\\
&=&\sum_{i=0}^{3}\left(q_i^T+\zeta_i+\frac{l_i-1}{2}\right)-\frac{l-1}{2}
=\sum_{i=1}^{3}\zeta_i
=\deg\phi_{\r}.
\end{eqnarray*}
\qed

As a consequence, we get a proof of Lemma \ref{deq}.

\noindent\emph{Proof of Lemma\ref{deq}.}
To begin with, Lemma \ref{GKZ-si} (with
$\zeta_i=0$) implies 
\beq\label{GKZ}
\prod_{i: \, l_i<0}\prod_{k=0}^{-l_i-1}\left(\delta+\frac{l}{l_i}
  (q_i^T+k)\right)\, \Phi
=C \si^l \prod_{i:\, l_i>0}\prod_{k=0}^{l_i-1}\left(\delta+\frac{l(q_i^T+k)}{l_i}\right)\,\Phi.
\eeq
The various values of $q_i^T$ and $l_i$ are listed in Tables \ref{PF-E6:weights}, \ref{PF-E7:weights} and 
\ref{PF-E8:weights}. 
We make the following observations:
\begin{enumerate}
\item If the RHS of \eqref{GKZ} contains a term $\delta+j$ with
  $j\in\mathbb{Z}, 1\leq j\leq l-1$, then the the reduced equation
  \eqref{PF-operator} has the form that we claimed.
\item For $l=3$, $\delta+1$ is always a factor of the RHS of \eqref{GKZ}.
\item If $l_i<0$, then for all $0\leq k\leq -l_i-1$,
  $l-\frac{l(q_i^T+k)}{l_i}$ is always a factor of the RHS of \eqref{GKZ}.
\end{enumerate}
This completes the proof of Lemma\ref{deq}.
\qed

The action of the operator \eqref{PF-operator} on a period integral is again a
period integral. The latter is holomorphic at $\si=0$; therefore, if it
is in the kernel of the Bessel operator \eqref{Bessel}, it must be a polynomial in
$\si$. But a non-zero period integral cannot
be a polynomial. In other words the period $\Phi_\r(\si)$ is a solution
to the Picard-Fuchs equation corresponding to the differential
operator \eqref{PF-operator}. In particular, we can check 
\begin{lemma}\label{PF-equations}
If $W$ is a Fermat simple elliptic singularity.
Let $x=C\si^l$; then either 
\beq\label{1st-order}
(1-x)\frac{\d}{\d x}\Phi_{\r}=(\deg\phi_{\r})\,\Phi_{\r},
\eeq
or $\Phi_\r$ saitisfies a hypergeometric equation 
\beq\label{2nd-order}
x(1-x)\frac{\d^2\Phi_{\r}}{\d
  x^2}+\left(\gamma_{\r}-\left(1+\alpha_{\r}+\beta_{\r}\right)x\right)\frac{\d\Phi_{\r}}{\d
  x}-(\alpha_{\r}\beta_{\r})\,\Phi_{\r} =0,
\eeq
where the weights $(\alpha_\r,\beta_\r,\gamma_\r)$ follows from \eqref{PF-operator} and satisfies
\beq\label{exp-id}
\alpha_{\r}+\beta_{\r}-\gamma_{\r}=\deg\phi_{\r}.
\eeq
Moreover, for $\r=\bf{0}$, $\Phi_{\r}$ satisfies \eqref{exp-id} for all invertible simple elliptic singularity $W$.
\end{lemma}
The first part of the lemma and the identity \eqref{exp-id} are
corollaries of Lemma \ref{identity}.
Unfortunately, we do not have a general combinatorial rule to determine which
indexes $(i',k')$ and $(i'',k'')$ should appear in
\eqref{PF-operator}. In other words, the second part of the Lemma is proved by straightforward 
computation, case by case. 

\begin{example}
For $W=X_1^6+X_2^3+X_3^3$, $\phi_{\m}=X_1^4X_2$, since $r_3=0$, we write $\r=(r_1,r_2)$ instead of $(r_1,r_2,r_3)$. The weights of the hypergeometric equations for $\Phi_{\r}$ are 
\ben
\left(\alpha_{\r},\beta_{\r},\gamma_{\r}\right)=\left(\frac{1+r_1}{12}, \frac{7+r_1}{12},\frac{2-r_2}{3}\right).
\een
\begin{table}[h]\center
\caption{Weights of periods for Fermat $E_8^{(1,1)}$}
\label{weight-E8}
\begin{tabular}{|c|c|c|c|c|c|c|c|c|c|c|c|c|c|c|}
\hline
$\phi_\r$&$X_1$&$X_2$&$X_1^2$&$X_1X_2$&$X_1^3$&$X_1^2X_2$&$X_1^4$&$X_1^3X_2$\\\hline
  $\alpha_\r,\beta_\r,\gamma_\r$
  &$\frac{1}{6},\frac{2}{3},\frac{2}{3}$
  &$\frac{1}{12},\frac{7}{12},\frac{1}{3}$
  &$\frac{1}{4},\frac{3}{4},\frac{2}{3}$
  &$\frac{1}{6},\frac{2}{3},\frac{1}{3}$
  &$\frac{1}{3},\frac{5}{6},\frac{2}{3}$
  &$\frac{1}{4},\frac{3}{4},\frac{1}{3}$
  &$\frac{5}{12},\frac{11}{12},\frac{2}{3}$
  &$\frac{1}{3},\frac{5}{6},\frac{1}{3}$\\
\hline
\end{tabular}
\end{table}
\end{example}

\subsection{Givental's theory}\label{sec:givental}
The goal here is to define the total ancestor
potential $\mathscr{A}^{\rm SG}_{W}(\s)$ at semisimple point $\s$ for $W$. 
Following Givental, we introduce the vector space
$\H=\mathscr{Q}_W((z))$ of formal Laurent series in $z^{-1}$ with
coefficients in $\mathscr{Q}_W$, equipped with the symplectic structure
\ben
\Omega(f(z),g(z))= {\rm res}_{z=0} (f(-z),g(z))dz.
\een
Using the polarization $\H=\H_+\oplus \H_-$, where $\H_+=\mathscr{Q}_W[z]$ and
$\H_-=\mathscr{Q}_W[[z^{-1}]]z^{-1}$ we identify $\H$ with the cotangent bundle
$T^*\H_+$.

Let $s\in \S$ be a semi-simple point, i.e., the critical values $u_i$
of $F$ ($1\leq i\leq \mu$) form locally near $s$ a coordinate system. Let us
also fix a path from $0\in S$ to $s$, so that we have a fixed branch
of the flat coordinates. Then we have an isomorphism 
\ben
\Psi_\s:\C^\mu\to H,\quad e_i\mapsto \sqrt{\Delta_i}\,\d_{u_i} 
\een
where $\Delta_i$ is determined by $(\d/\d u_i,\d/\d
u_j)=\delta_{ij}/\Delta_i$. It is well known that $\Psi_\s$
diagonalizes the Frobenius multiplication and the residue pairing,
i.e., 
\ben
e_i\bullet e_j = \sqrt{\Delta_i}e_i\delta_{i,j},\quad (e_i,e_j)=\delta_{ij}.
\een
Let $\S_{\rm ss}$ be the set of all semi-simple points. The complement
$\mathcal{K}=\S\setminus{\S_{\rm ss}}$ is an analytic hypersurface
also known as the {\em caustic}. It corresponds to deformations, s.t.,
$F$ has at least one non-Morse critical point. By definition
\ben
\S_{\rm ss}\to {\operatorname{Hom}}_\C(\C^\mu,H),\quad \s\mapsto \Psi_\s
\een
is a multi-valued analytic map.

The system of differential equations (\ref{frob_eq1}) and
(\ref{frob_eq2}) admits a unique formal asymptotical solution of the type 
\ben
\Psi_\s R_\s(z) e^{U_\s/z},\quad R_\s(z)=1+R_{\s,1}z+R_{\s,2} z^2+\cdots
\een
where $U_\s={\rm diag}\left(u_1(\s),\dots,u_\mu(\s)\right)$ is a diagonal matrix  and $R_{\s,k} \in {\operatorname{Hom}}_\C(\C^\mu,\C^\mu)$.
We refer to \cite{Du}, \cite{G1} for more details and proofs.

\subsubsection{The quantization formalism}

Let us fix a Darboux coordinate system on $\H$ given by the linear
functions $q_k^i$, $p_{k,i}$ defined as follows: 
\ben
\f(z) = \sum_{k=0}^\infty \sum_{i\in\mathfrak{R}} \ (q_k^i\, \d t_i\, z^k + p_{k,i}\,dt_i\,(-z)^{-k-1})\quad \in \quad \H,
\een
where $\{dt_{i}\}_{i\in\mathfrak{R}}$ is a basis of $H$ dual to
$\{\d t_{i}\}$ with respect to the residue pairing. 

If ${R}=e^{A(z)}$, where $A(z)$ is an infinitesimal symplectic
transformation, then we define $\widehat{R}$ as follows. Since $A(z)$
is infinitesimal symplectic, the map $\f\in \H \mapsto A\f\in \H$
defines a Hamiltonian vector field with Hamiltonian given by the
quadratic function $h_A(\f) = \frac{1}{2}\Omega(A\f,\f)$. By
definition, the quantization of $e^A$ is given by the differential
operator $e^{\widehat{h}_A},$ where the quadratic Hamiltonians are
quantized according to the following rules: 
\ben
(p_{k',e'}p_{k'',e''})\sphat = \hbar\frac{\d^2}{\d q_{k'}^{e'}\d q_{k''}^{e''}},\quad
(p_{k',e'}q_{k''}^{e''})\sphat = (q_{k''}^{e''}p_{k',e'})\sphat = q_{k''}^{e''}\frac{\d}{\d q_{k'}^{e'}},\quad
(q_{k'}^{e'}q_{k''}^{e''})\sphat =q_{k'}^{e'}q_{k''}^{e''}/\hbar.
\een
Note that the quantization defines a projective representation of the Poisson Lie algebra of quadratic Hamiltonians:
\ben
[\widehat{F},\widehat{G}]=\{F,G\}\sphat + C(F,G),
\een
where $F$ and $G$ are quadratic Hamiltonians and the values of the cocycle $C$ on a pair of Darboux monomials is non-zero only in the following cases:
\beq\label{cocycle}
C(p_{k',e'}p_{k'',e''},q_{k'}^{e'}q_{k''}^{e''})=
\begin{cases}
1 & \mbox{ if } (k',e')\neq (k'',e''),\\
2 & \mbox{ if } (k',e')=(k'',e'').
\end{cases}
\eeq

\subsubsection{The total ancestor potential}
By definition, the Kontsevich-Witten tau-function is the following generating series:
\beq\label{D:pt}
\D_{\rm pt}(\hbar;q(z))=\exp\Big( \sum_{g,n}\frac{1}{n!}\hbar^{g-1}\int_{\overline{\M}_{g,n}}\prod_{i=1}^n (q(\psi_i)+\psi_i)\Big),
\eeq
where $q(z)=\sum_k q_k z^k,$ $(q_0,q_1,\ldots)$ are formal variables,
$\psi_i$ ($1\leq i\leq n$) are the first Chern classes of the
cotangent line bundles on $\overline{\M}_{g,n}.$ The function is
interpreted as a formal series in $q_0,q_1+{\bf 1},q_2,\ldots$ whose coefficients are Laurent series in $\hbar$..

Let $\s\in \mathcal{S}_{\rm ss}$ be {\em a semi-simple}
point. Motivated by Gromov--Witten theory Givental introduced the
notion of the {\em total ancestor potential} of a semi-simple
Frobenius structure (see \cite{G1}, \cite{G2}). In our settings, the definition takes the form
\beq\label{SG-ancestor}
\mathscr{A}^{\rm SG}_{W}(\s)=\mathscr{A}_\s(\hbar;\q) :=
\widehat{\Psi_\s}\,\widehat{R_\s}\,e^{\widehat{U_\s/z}}\,
\prod_{i=1}^\mu \D_{\rm pt}(\hbar\,\Delta_i(\s); \leftexp{\it i}{\bf q} (z) \sqrt{\Delta_i(\s)})
\eeq
where
\ben
\q(z)=\sum_{k=0}^\infty
\sum_{j\in\mathfrak{R}} q_k^j\, z^k\d t_j, \quad 
\leftexp{\it i}{\bf q} (z) = \sum_{k=0}^\infty \leftexp{\it i}{q}_k \,
z^k.
\een
The quantization $\widehat{\Psi_s}$ is interpreted as the change of variables
\beq\label{change}
\sum_{i=1}^\mu \leftexp{\it i}{\bf q}(z)e_i=\Psi_\s^{-1}\q(z)\quad \mbox{i.e.}\quad
\leftexp{\it i}{q}_k\sqrt{\Delta_i} =\sum_{j\in\mathfrak{R}} (\d u^i/\d t_j)\, q_k^j.
\eeq

\section{Geometric limits: GW theory and FJRW theory}\label{sec:3}

The proof of Theorem \ref{t1} relies on a certain reconstruction
property of the mirror GW invariants and FJRW invariants. All the invariants are defined by intersections of cohomologies on $\overline{\M}_{g,n}$ associated with some
Cohomological Field Theory (CohFT). According to Krawitz--Shen \cite{KS}, starting with a
certain initial set of 3- and 4-point genus-0 correlators, we can
determine the remaining invariants using only the axioms of a
CohFT. The goal in this section is to introduce the CohFTs relevant
for Theorem \ref{t1} and to compute explicitly the initial
data of correlators needed for the reconstruction. 

\subsection{Cohomological Field Theories}
Let $H$ be a vector space of dimension $N$ with a unit {\bf 1} and a non-degenerate paring $\eta$. Without loss of generality, we always fix a basis of $H$, say $\mathscr{S}:=\{\d_i, i=0,\dots,N-1\}$, and we set $\d_0={\bf 1}$. Let $\{\d^j\}$ be the dual basis in the dual space $H^{\vee}$,
(i.e., $\eta(\d_i,\d^j)=\delta_i^j$).
A CohFT $\Lambda$ is a set of multi linear maps $\{\La_{g,n}\}$, for each stable genus $g$ curve with $n$ marked points, i.e., $2g-2+n>0$,
$$\Lambda_{g,n}: H^{\otimes n}\longrightarrow H^*(\overline{\M}_{g,n}, \C).$$

Furthermore, $\La$ satisfies a set of axioms (CohFT axioms) described below:
\begin{enumerate}
\item ($S_n$-invariance) For any $\sigma\in S_n$, and $\gamma_1,\dots,\gamma_n\in H$, then
$$\Lambda_{g,n}(\gamma_{\sigma(1)}, \dots, \gamma_{\sigma(n)})=\La_{g,n}(\gamma_1,\dots,\gamma_n).$$
\item (Gluing tree)
Let
$$\rho_{tree}:\overline{\M}_{g_1,n_1+1}\times\overline{\M}_{g_2,n_2+1}\to\overline{\M}_{g,n}$$
where $g=g_1+g_2, n=n_1+n_2$, be the morphism induced from gluing the last marked point of the first curve and the first marked point of the second curve; then
\begin{align*}
\rho_{tree}^*&\big(\Lambda_{g,n}(\gamma_1,\dots,\gamma_n)\big)\\
&=\sum_{\alpha,\beta\in\mathscr{S}}\Lambda_{g_1,n_1+1}(\gamma_1,\dots,\gamma_{n_1},\alpha)\eta^{\alpha,\beta}\Lambda_{g_2,n_2+1}(\beta,\gamma_{n_1+1},\dots,\gamma_n).
\end{align*}
Here $\big(\eta^{\alpha,\beta}\big)_{N\times N}$ is the inverse matrix of $\big(\eta(\alpha,\beta)\big)_{N\times N}$.
\item (Gluing loop)
Let
$$\rho_{loop}:\overline{\M}_{g-1,n+2}\to\overline{\M}_{g,n},$$
be the morphism induced from gluing the last two marked points; then
\beqa
\rho_{loop}^*\big(\Lambda_{g,n}(\gamma_1,\dots,\gamma_n)\big)=\sum_{\alpha,\beta\in\mathscr{S}}\Lambda_{g-1,n+2}(\gamma_1,\dots,\gamma_n,\alpha,\beta)\eta^{\alpha,\beta}.
\eeqa
\item (Pairing)
$$\int_{\overline{\M}_{0,3}}\Lambda_{0,3}({\bf 1}, \gamma_1, \gamma_2)=\eta(\gamma_1, \gamma_2).$$
\end{enumerate}
If in addition the following axiom holds:
\begin{itemize}
\item[(5)] (Flat identity) Let
  $\pi:\overline{\M}_{g,n+1}\to\overline{\M}_{g,n}$ be the forgetful morphism; then
$$\Lambda_{g, n+1}(\gamma_1,\dots,\gamma_n, {\bf 1})=\pi^* \Lambda_{g,n}(\gamma_1,\dots,\gamma_n).$$
\end{itemize}
then we say that $\La$ is a CohFT with a {\em flat identity}.

If $\Lambda$ is a CohFT; then there is a natural formal family of
CohFTs $\Lambda(\t)$. Namely, 
\ben
\Lambda(\t)_{g,n}(\gamma_1,\dots,\gamma_n) = \sum_{k=0}^\infty
\frac{1}{k!}\, \pi_*\left(\Lambda_{g,n+k}(\gamma_1,\dots,\gamma_n,\t,\dots,\t)\right),
\een
where $\pi:\overline{\M}_{g,n+k}\to \overline{\M}_{g,n}$ is the
morphism forgetting the last $k$ marked points. 
Note that $\La(\t)_{0,3}$ induces a family of Frobenius multiplications $\bullet_\t$ on $(H,\eta)$, defined by
\beq\label{frobenius}
\eta(\alpha\bullet_\t\beta,\gamma)=\int_{\overline{\M}_{0,3}}\Lambda(\t)_{0,3}(\alpha,\beta,\gamma).
\eeq
It is well known that \cite{M} the genus-$0$ part of the CohFT $\{\Lambda(\t)_{0,n}\}$ is equivalent to a Frobenius
manifold $(H,\eta,\bullet_\t)$, in the sense of Dubrovin \cite{Du}. We call the vector space $H$ the {\em state space} of the
CohFT. 

\subsubsection{The total ancestor potential of a CohFT}
For a given CohFT $\Lambda$, the {\em correlator functions} are by definition the
following formal series in $\t\in H$:
\beq\label{cor-anc}
\LD
\tau_{k_1}(\alpha_1),\dots,\tau_{k_n}(\alpha_n)\RD_{g,n}^{\Lambda}(\t)=
\int_{\overline{\M}_{g,n}}
\Lambda(\t)_{g,n}(\alpha_1,\dots
,\alpha_n)\,\psi_1^{k_1}\dots\psi_n^{k_n}\, ,
\eeq
where $\psi_i$ is the $i$-th psi class on ${\overline{\M}_{g,n}}$, $\alpha_i\in H,$ and $k_i\in \Z_{\geq 0}$. 
The value of a correlator function at $\t=0$ is called simply a
{\em correlator}. We call $g$ the genus of the correlator
function and each $\tau_{k_i}(\alpha_i)$ is called a {\em
descendant} (resp. {\em non-descendant}) {\em insertion} if $k_i>0$ (resp. $k_i=0$). 

For each basis $\{\d_i\}$ in $H$, we fix a sequence of formal
variables $\{q_k^i\}_{k=0}^\infty$ and define
\ben
\q(z) = \sum_{k=0}^\infty \sum_{i=0}^{N-1}\ q_k^i\, \d_{i}\, z^k\ ;
\een
then the \emph{genus-$g$ ancestor potential} is the following
generating function:
\begin{equation*}
\mathscr{F}_g^{\Lambda(\t)}(\q):=
\sum_{n}\frac{1}{n!}
\LD\q(\psi_1)+\psi_1,\dots,\q(\psi_n)+\psi_n\RD_{g,n}^{\Lambda}(\t),
\end{equation*}
where each correlator should be expanded multi linearly in $\q$ and
the resulting correlators are evaluated according to
\eqref{cor-anc}. Let us point out that we have assumed that the CohFT has a
flat identity ${\bf 1}\in H$ and we have incorporated the {\em
  dilaton shift} in our function, so that $\mathscr{F}_g^{\Lambda(\t)}$ is a
formal series in $q_k$, $k\neq 0$ and $q_1+{\bf 1}$.  
Finally, the \emph{total ancestor potential of a CohFT $\Lambda(\t)$} is defined by
\begin{equation}\label{ancestor}
\mathscr{A}^{\Lambda(\t)}\left(\hbar;\q\right):=\exp\left(\sum_{g=0}^{\infty}\hbar^{2g-2}\mathscr{F}_g^{\Lambda(\t)}(\q)\right).
\end{equation}

\subsection{GW theory of elliptic orbifold $\mathbb{P}^1$}\label{elliptic-GW}
Let $\mathcal{X}:=\mathbb{P}^1_{a_1,a_2,a_3}$ be the orbifold-$\mathbb{P}^1$ with 3
orbifold points, such that, the i-th one has isotropy group
$\mathbb{Z}/a_i\mathbb{Z}$. We are interested in 3 cases:
$(a_1,a_2,a_3)$$=(3,3,3),(4,4,2),(6,3,2)$.  Together with $\mathbb{P}^1_{2,2,2,2}$, they correspond to
orbifold-$\mathbb{P}^1$s that are quotients of an elliptic curve by a
finite group. 
The Chen-Ruan cohomology $H^*_{\rm CR}(\mathbb{P}^1_{a_1,a_2,a_3})$ has
the following form:
\beqa
H^*_{\rm CR}(\mathbb{P}^1_{a_1,a_2,a_3})=\bigoplus_{i=1}^3\bigoplus_{j=1}^{a_i-1}\mathbb{C}[\Delta_{ij}]\bigoplus\mathbb{C}[\Delta_{01}]\bigoplus\mathbb{C}[\Delta_{02}].
\eeqa
where $\Delta_{01}=1$ and $\Delta_{02}$ is the Poincar\'e dual to a
point. The classes $\Delta_{ij}(1\leq i\leq 3,1\leq j\leq a_i-1)$ are
in one-to-one correspondence with the twisted sectors. The latter are
just orbifold points, and we define $\Delta_{ij}$ to be the unit in the cohomology
of the corresponding twisted sector. Our indexing matches the
complex degrees 
\ben
\deg\Delta_{ij}={j}/{a_i},\quad 1\leq i\leq 3, \quad 1\leq j\leq a_i-1.
\een
$\mathcal{X}$ is a compact Kahl\"er orbifold. Let
$\overline{\mathcal{M}}_{g,n,d}^{\mathcal{X}}$ be the moduli space
of degree-$d$ stable maps from a genus-$g$ orbi-curve, equipped with
$n$ marked points, to $\mathcal{X}$. Let us denote by 
$\pi$ the forgetful map, and by ${\rm
  ev}_i$ the evaluation at the $i$-th marked point
\ben
\begin{array}[c]{ccccc}
\overline{\mathcal{M}}_{g,n}&
\stackrel{\pi}{\longleftarrow}&
\overline{\mathcal{M}}_{g,n+k,d}^{\mathcal{X}}&
\stackrel{\rm ev_i}{\longrightarrow}&
I\mathcal{X}\ .
\end{array}
\een
The moduli space is
equipped with a virtual fundamental cycle
$[\overline{\mathcal{M}}_{g,n,d}^{\mathcal{X}}]$, such that the maps
$(\Lambda^\mathcal{X})_{g,n}:
H^*_{\rm CR}(\mathbb{P}^1_{a_1,a_2,a_3})^{\otimes n} \longrightarrow
H^*(\overline{\M}_{g,n};\C)$
defined by 
\ben
(\Lambda^\mathcal{X})_{g,n}(\alpha_1,\dots,\alpha_n):=
\pi_*\, \Big(
[\overline{\mathcal{M}}_{g,n,d}^{\mathcal{X}}]  \cap 
\prod_{i=1}^n {\rm ev}_i^*(\alpha_i)\ 
\Big)
\een
form a CohFT with state space
$H^*_{CR}(\mathbb{P}^1_{a_1,a_2,a_3})$. The total ancestor potential
$\mathscr{A}^{\rm GW}_\mathcal{X}$  of $\mathcal{X}$ is by definition \eqref{ancestor} the total ancestor
potential of the CohFT $\Lambda^{\mathcal{X}}(\t={\bf 0})$. For more details on orbifold Gromov--Witten theory, we refer to \cite{CheR}.
For $\mathbb{P}^1_{a_1,a_2,a_3}$, the orbifold Poincar\'e pairing takes the form 
\beq\label{CR-pair}
\LD\Delta_{i_1j_1},\Delta_{i_2j_2}\RD=\left\{\begin{array}{ll}
\left(\delta_{i_1,i_2}\delta_{j_1+j_2,a_k}\right)/a_k, & k=i_1, i_1+i_2\neq0;\\
\delta_{j_1+j_2,3},&i_1=i_2=0.
\end{array}\right.
\eeq
It is easy to compute that the above 3-point correlators are 
\beq\label{CR-product}
\LD\Delta_{i_1j_1},\Delta_{i_2j_2},\Delta_{i_3j_3}\RD_{0,3,0}=
\left\{
\begin{array}{ll}
1/a_k,
&
i_1=i_2=i_3=k, j_1+j_2+j_3=a_k;
\\
\LD\Delta_{i_2j_2},\Delta_{i_3j_3}\RD,
& 
(i_1,j_1)=(0,1);
\\
0,
&
\text{otherwise}.
\end{array}
\right.
\eeq
Recall the Chen-Ruan orbifold cup product $\bullet$ on
$H^*_{CR} ( \mathbb{P}^1_{a_1,a_2,a_3};\C)$ is defined by pairing and 3-point correlators \eqref{frobenius}.
According to Krawitz and Shen \cite{KS}, we have the following reconstruction result. 
\begin{lemma}\label{thm:GW recons}
The Gromov-Witten total ancestor potential $\mathscr{A}^{\rm GW}_\mathcal{X}$ of elliptic orbifold  $\mathcal{X}=\mathbb{P}^1_{3,3,3}$,
$\mathbb{P}^1_{4,4,2},\mathbb{P}^1_{6,3,2}$ is determined by the following initial data: the Poincar\'e pairing, the Chen-Ruan product, and the correlator 
\beq\label{3-pt-deg1}
\langle\Delta_{1,1},\Delta_{2,1},\Delta_{3,1}\rangle_{0,3,1}^{\mathcal{X}}=1.
\eeq
\end{lemma}
In particular, using Lemma \ref{thm:GW recons}, one can construct a
mirror map to identify the Gromov--Witten theory of
the above orbifolds and the Saito--Givental theory associated to
certain ISES (see \cite{KS}, \cite{MR}). The genus 0 reconstruction of the GW theory for those orbifolds are obtained by Satakea and Takahashi independently \cite{ST}. Moreover, they also proved the isomorphism between the Frobenius manifold from Gromov-Witten theory of $\mathbb{P}^1_{3,3,3}$ and Saito's Frobenius manifold of $X_1^3+X_2^3+X_3^3+\sigma X_1X_2X_3$ at $\sigma=\infty$.

On the other hand, since the mirror
symmetry identifies the correlation functions with certain period
integrals, by analyzing the monodromy of the period integrals, one can
prove that the Gromov-Witten invariants are quasi-modular forms on
some finite index subgroups of the modular group, \cite{MR}, \cite{MS}. 
In particular,
the non-zero, genus-0, 3-point
correlators are modular forms of weight 1. Let us list
the first few terms of their Fourier series. For $\mathcal{X}=\P^1_{3,3,3}$, the following correlators are weight-1 modular forms
on $\Gamma(3)$,
\ben
\left\{
\begin{aligned}
\LD\Delta_{1,1},\Delta_{2,1},\Delta_{3,1}\RD_{0,3}&=q+q^4+2q^7+2q^{13}+\cdots=\frac{\eta(9\tau)^3}{\eta(3\tau)}\\
\LD\Delta_{1,1},\Delta_{1,1},\Delta_{1,1}\RD_{0,3}&=1/3+2q^3+2q^9+2q^{12}+\cdots=\frac{3\eta(9\tau)^3+\eta(\tau)^3}{3\eta(3\tau)},
\end{aligned}
\right.
\een
For $\mathcal{X}=\P^1_{4,4,2}$, the following correlators are weight-1 modular forms on
$\Gamma(4)$,
\ben
\left\{
\begin{aligned}
\LD\Delta_{1,1},\Delta_{2,1},\Delta_{3,1}\RD_{0,3}&=q+2q^5+q^9+2q^{13}+\cdots\\
\LD\Delta_{1,1},\Delta_{1,1},\Delta_{1,2}\RD_{0,3}&=1/4+q^4+q^8+q^{16}+\cdots\\
\LD\Delta_{1,1},\Delta_{1,1},\Delta_{2,2}\RD_{0,3}&=q^2+2q^{10}+q^{18}+\cdots\ ,
\end{aligned}
\right.
\een
For $\mathcal{X}=\P^1_{6,3,2}$, the following correlators are weight-1 modular forms on
$\Gamma(6)$ 
\ben
\left\{
\begin{aligned}
\LD\Delta_{1,1},\Delta_{2,1},\Delta_{3,1}\RD_{0,3}&=q+2q^7+2q^{13}+2q^{19}+\cdots\\
\LD\Delta_{1,1},\Delta_{1,1},\Delta_{1,4}\RD_{0,3}&=1/6+q^6+q^{18}+q^{24}+\cdots\\
\LD\Delta_{1,1},\Delta_{1,2},\Delta_{1,3}\RD_{0,3}&=1/6+q^{12}+\cdots\\
\LD\Delta_{1,1},\Delta_{1,1},\Delta_{2,2}\RD_{0,3}&=q^2+q^8+2q^{14}+\cdots\\
\LD\Delta_{1,1},\Delta_{1,2},\Delta_{3,1}\RD_{0,3}&=q^3+q^9+2q^{21}+\cdots\\
\LD\Delta_{1,1},\Delta_{2,1},\Delta_{1,3}\RD_{0,3}&=q^4+q^{16}+2q^{28}+\cdots\\
\end{aligned}
\right.
\een
\begin{remark}
The $q$-expansions above are mainly obtained by WDVV equations (see explicit formulas in the Appendix of \cite{S}).
In each of the above cases the given modular forms form a basis in the space of modular forms of
weight 1. One can prove that these modular forms give an embedding of
the corresponding modular curve in a projective space.
\end{remark}

\subsection{FJRW theory of invertible simple elliptic singularities}\label{sec-3.2}
For any non-degene\-rate, quasi-homogeneous polynomial $W$,
Fan--Jarvis--Ruan, following a suggestion of Witten, introduced a
family of moduli spaces and constructed a virtual fundamental cycle.
The latter gives rise to a cohomological field theory, which is now
called the FJRW theory. 
In our paper, we make use of the FJRW theories associated with the
invertible simple elliptic singularities, i.e., $W$ is one of the
polynomials in Table 1. Let us briefly  review the FJRW theory only for such
$W$ and refer to \cite{FJR} for the general case and more details.

\subsubsection{FJRW vector space and axioms}
Recall the \emph{group of diagonal symmetries} $G_{W}$ of the polynomial $W$ is
\ben
G_{W}:=
\left\{(\la_1,\la_2,\la_3)\in\left(\mathbb{C}^*\right)^3\Big\vert\,W(\la_1\,X_1,\la_2\,X_2,\la_3\,X_3)=W(X_1,X_2,X_3)\right\}.
\een
The \emph{FJRW state space} $\mathscr{H}_{W,G_{W}}$ (or $\mathscr{H}_{W}$ for short) is the direct sum of all $G_{W}$-invariant relative cohomology:
\begin{equation}
\mathscr{H}_{W}:=
\bigoplus_{h\in G_{W}}H_{h},\quad
H_h:=H^*(\mathbb{C}_{h};W^{\infty}_{h};\mathbb{C})^{G_{W}}.
\end{equation}
Here $\mathbb{C}_{h}(h\in G_{W})$ is the $h$-invariant subspace of
$\mathbb{C}^3$, $W_{h}$ is the restriction of $W$ to
$\mathbb{C}_{h}$ and ${\rm Re}W_{h}$ is the real part of $W_{h}$, and $W^{\infty}_{h}=(\mathrm{Re}W_{h})^{-1}(M,\infty)$, for some $M\gg0.$

The vector space $H_h(h\in G_{W})$ has a natural grading given by
the degree of the relative cohomology classes. However, for the
purposes of the FJRW theory we need a modification of the standard
grading. Namely, for any $\alpha\in H_{h}$, we define (this is half of
$\deg_W\alpha$ in \cite{FJR}) 
$$
\deg_{W}\alpha:=\frac{\deg\alpha}{2}+\sum_{i=1}^3(\Theta_i^{h}-q_i),
$$
where ${\rm deg}\ \alpha$ is the degree of
$\alpha$ as a relative cohomology class in $H^*(\mathbb{C}_{h};W^{\infty}_{h};\mathbb{C})$
and 
\ben
h=\left(\be[\Theta_1^{h}],\be[\Theta_2^{h}],\be[\Theta_3^{h}]\right)\in(\mathbb{C}^*)^3, \quad \Theta_i^{h}\in[0,1)\cap\mathbb{Q}
\een 
where for $y\in\mathbb{R}$, we put
$\be[y]:=\exp(2\pi\sqrt{-1}y)$. Clearly the numbers $\Theta_i^h$ are uniquely
determined from $h$. For any $\alpha\in H_h$, we define
\begin{equation}\label{theta-weight}
\Theta(\alpha):=\left(\be[\Theta_1^{h}],\be[\Theta_2^{h}],\be[\Theta_3^{h}]\right).
\end{equation}
The elements in $H_h$ are called \emph{narrow} (resp. \emph{broad})
and $H_h$ is called a {\em narrow sector} (resp. {\em broad sector}) if
$\C_h=\{0\}$ (resp. $\C_h\neq\{0\}$). For invertible simple elliptic
singularities, the space $H^*(\C_h;W_{h}^{\infty};\mathbb{Q})$ is one-dimensional for all narrow sectors $H_h$. We always choose a generator
$\alpha\in H_{h}$ such that
\beq\label{narrow}
\alpha:=1\in H^*(\C_h;W_{h}^{\infty};\mathbb{Q})^{G_W}.
\eeq
In general, in order to describe the broad sectors, we have to represent the relative cohomology classes by
differential forms; then there is an identification (see \cite{FJR} and the references there)
\beq\label{FJRW-pairing}
\Big( \mathscr{H}_{W,G},\langle\ ,\ \rangle\Big)\equiv\Big( \bigoplus_{h\in
    G}\left(\mathscr{Q}_{W_h}\omega_h\right)^G,{\rm Res}\Big)\, ,
\eeq
where $\omega_h$ is the restriction of the standard volume form to the fixed
locus $\C_h$, ${\rm Res}$ is the residue pairing, and 
$\langle\ ,\ \rangle$ is a non-degenerate pairing induced from the
intersection of relative homology cycles.
There exists a basis of the narrow sectors such that the
pairing $\langle v_1 , v_2 \rangle$, $v_i\in H_{h_i}$, is $1$ if
$h_1h_2=1$ and $0$ otherwise. The vectors in the broad sectors are
orthogonal to the vectors in the narrow sectors. In order to compute
the pairing on the broad sectors one needs to
use the identification \eqref{FJRW-pairing} and compute an appropriate
residue pairing. In our case however, we can
express all invariants using narrow sectors only. So a more detailed
description of the broad sectors is not needed. We refer to \cite{FJR}
for more details.

Let $(W,G)$ be an admissible pair. A {\em $W$-spin structure} on a genus-$g$
Riemann surface $C$ with $n$ marked orbifold points $(z_1,\dots,z_n)$ is a
collection of $N$ ($N$ is the number of variables in $W$) orbifold
line bundles $\mathcal{L}_1, \dots,
\mathcal{L}_N$ on $C$ and isomorphisms 
\ben
\psi_a:M_a(\mathcal{L}_1, \dots, \mathcal{L}_N)\to
\omega_C(-z_1-\cdots-z_n),
\een
where $\omega_C$ is the dualizing sheaf on $C$ and $M_a$ are the
homogeneous monomials whose sum is $W$. The orbifold line bundles have
a monodromy near each marked point $z_i$ which determines an element
$h_i\in G$. In particular, if $H_{h_i}$ is a narrow (resp. broad) sector we say
that the marked point is narrow (resp. broad).  For fixed $g,n,$ and $h_1,\dots,h_n\in G$, 
Fan-Jarvis-Ruan (see \cite{FJR}) constructed the compact moduli space $\overline{\mathscr{W}}_{g,n}(h_1,\cdots,h_n)$ 
of nodal Riemann surfaces equipped with a $W$-spin structure. 
In this compactification the line bundles $(\L_1,\dots,\L_N)$ are allowed to be orbifold at
the nodes in such a way that the monodromy around each node is an element
of $G$ as well.
The moduli space has a decomposition into a disjoint union of moduli
subspaces $\overline{\mathscr{W}}_{g,n}(\Gamma_{h_1,\dots,h_n})$ consisting of
$W$-spin structures on curves $C$ whose dual graph is
$\Gamma_{h_1,\dots,h_n}$. 
Recall that the dual graph of a nodal curve $C$ is a graph whose vertices are the
irreducible components of $C$, edges are the nodes, and tails are the
marked points. The latter are decorated by elements $h_i\in G$, so the
tails of our graphs are also colored respectively. 
We omit the
subscript $(h_1,\dots,h_n)$ whenever the decoration is understood from
the context. The connected component
$\overline{\mathscr{W}}_{g,n}(\Gamma_{h_1,\dots,h_n})$ is naturally stratified
by fixing the monodromy
transformations around the nodes, i.e., the strata are in one-to-one
correspondence with the colorings of the edges of the dual graph $\Gamma_{h_1,\cdots,h_n}$.

Fan--Jarvis--Ruan constructed a virtual fundamental cycle
$[\overline{\mathscr{W}}_{g,n}(\Gamma)]^{\mathrm{vir}}$ of
$\overline{\mathscr{W}}_{g,n}(\Gamma)$ (see \cite{FJR1}),
which gives rise to a CohFT
$$\Lambda_{g,n}^{W,G}:\left(\mathscr{H}_{W,G}\right)^{\otimes n}\longrightarrow H^*(\overline{\mathcal{M}}_{g,n}).$$
For brevity put $\Lambda_{g,n}^{W}$ for $\Lambda_{g,n}^{W,G_{W}}$. We define \emph{the total ancestor FJRW potential} $\mathscr{A}_{W}^{\rm FJRW}$ of $(W,G_W)$ from the CohFT $\Lambda^{W}$ using \eqref{ancestor}.  

Finally, let us list some general properties of the FJRW correlators of a simple
elliptic singularity $W$, see \cite{FJR} for the proofs.
\begin{enumerate}
\item (Selection rule)
If the correlator
$\LD\tau_{k_1}(\alpha_1),\dots,\tau_{k_n}(\alpha_n)\RD_{g,n}^{W}$
is non-zero; then 
\begin{equation}\label{deg:FJRW}
\sum_{i=1}^n\deg_{W}(\alpha_i)+\sum_{i=1}^n k_i=2g-2+n.
\end{equation}

\item (Line bundle criterion). If the moduli space $\overline{\mathscr{W}}_{g,n}(h_1,\dots,h_n)$ is non-empty, then the degree of the desingularized line bundle $|\mathcal{L}_j|$ is an integer, i.e.
\begin{equation}\label{deg:Line BD}
\deg(|\mathcal{L}_j|)=q_j(2g-2+n)-\sum_{k=1}^{n}\Theta_{j}^{h_k}\in\mathbb{Z}.
\end{equation}

\item (Concavity)
Suppose that all marked points are narrow, $\pi$ is the morphism from the universal curve to $\overline{\mathscr{W}}_{g,n}(h_1,\dots,h_n)$ and $\pi_{*}\bigoplus_{i=1}^{3}\mathcal{L}_i=0$ holds; then
\begin{equation}\label{eq:Concave}
[\overline{\mathscr{W}}_{g,n}(h_1,\dots,h_n)]^{\mathrm{vir}}
=c_{\textrm{top}}\left(-R^1\pi_*\bigoplus_{i=1}^3\mathcal{L}_i\right)\cap[\overline{\mathscr{W}}_{g,n}(h_1,\dots,h_n)].
\end{equation}

\end{enumerate}

Let $\alpha_i=1\in H_{h_i}$, $1\leq i\leq 4$ be the generators
(cf. \eqref{narrow}). The concavity formula \eqref{eq:Concave}
implies that $\Lambda_{0,4}^{W}(\alpha_1,\dots,\alpha_4)\in
H^*(\overline{\M}_{0,4},\mathbb{C})$. According to the orbifold
Grothendieck-Riemann-Roch formula 
(see \cite{C},  Theorem 1.1.1),
\begin{equation}\label{eq:OGRR}
\Lambda_{0,4}^{W}\left(\alpha_1,\dots,\alpha_4\right)=\sum_{i=1}^3\left(\frac{B_{2}(q_i)}{2}\kappa_1
-\sum_{j=1}^4\frac{B_{2}(\Theta_{i}^{h_j})}{2}\psi_j
+\sum_{\Gamma\in\mathbf{\Gamma}_{0,4,W}(h_1,\dots,h_4)}
\frac{B_{2}(\Theta_{i}^{h_{\Gamma}})}{2}[\Gamma]\right),
\end{equation}
where  $B_2$ is the second Bernoulli polynomial
$$B_2(y)=y^2-y+\frac{1}{6},$$
$[\Gamma]$ is the boundary class on $\overline{\M}_{g,n}$ corresponding
to the graph $\Gamma$, and
$\mathbf{\Gamma}_{0,4,W}(h_1,\dots,h_4)$ is the set of
graphs with one edge decorated by $G_{W^T}$. The graph $\Gamma$ has 4
tails decorated by $h_1,h_2,h_3,h_4$ and its edge is decorated by
$h_{\Gamma}$ and $h_{\Gamma}^{-1}$. 
If the moduli space $\overline{\mathscr{W}}_{0,4}(h_1,\dots,h_4)$ is non-empty, each component satisfies
\eqref{deg:Line BD}. It is easy to see that the formula does not
depend on the choice of $h_{\Gamma}$ or $h_{\Gamma}^{-1}$.

\subsubsection{Generators of the FJRW ring.}
From now on, we will consider $W$ as an ISES in Table \eqref{table-1}. 
Let $(W^T,G_{W^T})$ be the Berglund-H\"ubsch-Krawitz mirror of $(W,\{1\})$ (see its definition in Section \ref{gms-int}). We will compare
the FJRW theory for $(W^T,G_{W^T})$ with the Saito-Givental theory for $W$. 

Since $\mathscr{H}_{W^T}:=\mathscr{H}_{W^T,G_{W^T}}$ is the state space of a CohFT, it has a
Frobenius algebra structure, where the multiplication $\bullet$ is defined by pairing and 3-point correlators \eqref{frobenius}. 
For all invertible $W$,  M. Krawitz (see \cite{K}) constructed a ring isomorphism 
\beq\label{ring-iso}
\mathscr{H}_{W^T}\cong\mathscr{Q}_W.
\eeq 
Next, we give an explicit description of the generators of $\mathscr{H}_{W^T}$ and the ring isomorphism for $W$ is an ISES in Table \eqref{table-1}. 
For a more general description, we refer the interested readers to \cite{K}, \cite{FS} and \cite{A}.

For every ISES $W^T$, there exists a $3$-tuple $(a,b,c)\in\Z^3$ such that
$$
G_{W^T}\cong\mu_{a}\times\mu_{b}\times\mu_{c},\quad \mu_k=\Z/k\Z.
$$
We assume $a\geq b\geq c$ and omit the factor $\mu_1$. For example, for $W^T=X_1^3+X_1X_2^4+X_3^2$
$$
G_{W^T}=\left\{(\la_1,\la_2,\la_3)\Big\vert\la_1^3=\la_1\la_2^4=\la_2^2=1\right\}\cong\mu_{12}\times\mu_{2}.$$
Let $h=(i,j,k)\in\mu_{a}\times\mu_{b}\times\mu_{c}\cong G_{W^T}$, $H_{h}$ is a narrow sector and $H_h\cong\C$, if
$$1\leq i<a, \quad 1\leq j<b,\quad 1\leq k<c,$$
In this case, we denote a generator of $H_h$ by
$$\be_{i,j,k}:=1\in H_h=H^0(\C_h;W_h^{\infty};\mathbb{Q}).$$ 

\begin{example}\label{ex:loop}
We compute the FJRW ring for loop singularity $W^T$, with $W\in E_{6}^{(1,1)}$. 
$$W^T=X_1^2X_3+X_1X_2^2+X_2X_3^2,\quad G_{W^T}=\left\{\be_i=\left(\be[\frac{i}{9}],\be[\frac{4i}{9}],\be[-\frac{2i}{9}]\right),i=1,\dots,8\right\}\cong\mu_{8}.$$
All nonzero 3-point genus-0 correlators are
$$\left\{
\begin{array}{l}
\LD\be_1,\be_1,\be_1\RD_{0,3}=\LD\be_4,\be_4,\be_4\RD_{0,3}=\LD\be_7,\be_7,\be_7\RD_{0,3}=-2;\\
\LD\be_3,\be_i,\be_{9-i}\RD_{0,3}=\LD\be_1,\be_4,\be_7\RD_{0,3}=1.\\
\end{array}
\right.
$$
The first row uses Index Zero Axiom (see \cite{FJR}) and the second row uses Concavity Axiom \eqref{eq:Concave}. It is easy to see $\be_3$ is the identity element and the ring relations are
$$2\,\be_1\bullet\be_4+\be_7^2=2\,\be_4\bullet\be_7+\be_1^2=2\,\be_7\bullet\be_1+\be_4^2=0.$$
Thus we obtain a ring isomorphism between $\mathscr{H}_{W^T}$ and $\mathscr{Q}_{W}$:
$$\rho_1=\be_4\mapsto\,X_1,\quad \rho_2=\be_1\mapsto\,X_2,\quad\rho_3=\be_7\mapsto\,X_3.$$
\end{example}

For all 13 types of ISESs with a maximal admissible group, 
there is a unique narrow sector $\rho_{-1}$, with $\deg_{W^T}(\rho_{-1})=1$ and 
$$\Theta(\rho_{-1}):=\left(1-q_1^{T}, 1-q_2^T,1-q_3^T\right).$$

There are 13 types of ISESs, but only 9 of them do not have broad generators. The narrow sectors have the advantage that we can use the
powerful concavity axiom \eqref{eq:Concave}. Combined with the remaining
properties of the correlators and the WDVV
equations this allows us to reconstruct all genus-0 FJRW
invariants. According to the reconstruction theorem in \cite{KS}, we
can also reconstruct the higher genus FJRW invariants, i.e., the total ancestor potential function
$\mathscr{A}^{\rm FJRW}_{W^T}.$ 

In the remaining 4 cases, we know how to offset the complication of
having broad generators for $W^T=X_{1}^2+X_1X_2^2+X_2X_3^3.$ 
The maximal abelian group is of order $12$. Its FJRW vector space has eight generators:
$$\be_1,\be_3,\be_5,\be_7,\be_9,\be_{11},R_4,R_8.$$
Here $R_4$ and $R_8$ are the cohomology classes represented by the following forms: 
$$R_h=dX_1\wedge dX_2\in H^2(\C_{h};W_{h}^{\infty};\mathbb{Q}),\quad  h=4,8\in G_{W^T}.$$
Note that $R_4$ and $R_8$ are $G_{W^T}$-invariant elements in
$\mathscr{Q}_{W_h}\omega_h$ where $h\in G_{W^T}$ acts on each
factor $X_i$ and $dX_i$ as
multiplication by $\be[q_i^T]$. Although one of the ring generators
($R_4$) is broad, we have enough WDVV equations to reconstruct the correlators containing broad sectors from correlators with only narrow elements and apply the concavity axioms. 

For the other three types of ISESs, we can still compute some genus-0
4-point correlators with broad sectors, but we do not know how to reconstruct the
complete theory only from correlators with narrow elements. 
In other words, for 10 out of the 13 ISESs, we can compute all the FJRW
invariants. These cases and the corresponding ring generators $\rho_i$ of
the FJRW ring $\mathscr{H}_{W^T}$ are listed in tables below. 

\begin{table}[h]\center
\caption{Generators of the FJRW ring $\mathscr{H}_{W^T}, \quad W\in E^{(1,1)}_6$}
\label{ring-E6}
\begin{tabular}{|c|c|c|ccc|c|}
\hline
$W^T$ &$G_{W^T}$& $\Theta(\be_{i,j,k})$&$\rho_1$&$\rho_2$&$\rho_3$&$\rho_{-1}$  \\
\hline

$X_1^3+X_2^3+X_3^3$ &$\mu_{3}^3$& $\be[\frac{i}{3}],\be[\frac{j}{3}],\be[\frac{k}{3}]$&$\be_{2,1,1}$&$\be_{1,2,1}$&$\be_{1,1,2}$&$\be_{2,2,2}=\rho_1\rho_2\rho_3$  \\

\hline

$X_1^2X_3+X_1X_2^2+X_2X_3^2$ &$\mu_{8}$& $\be[\frac{i}{9}],\be[\frac{4i}{9}],\be[-\frac{2i}{9}]$&$\be_{4}$&$\be_{1}$&$\be_{7}$&$\be_{6}=\rho_1\rho_2\rho_3$  \\
\hline

$X_1^2+X_1X_2^2+X_2X_3^3$ &$\mu_{12}$& $\be[\frac{i}{2}],\be[-\frac{i}{4}],\be[\frac{i}{12}]$&$R_4$&$\be_{1}$&$\be_{7}$&$\be_{9}=\rho_2\rho_3^2$  \\

\hline

\end{tabular}
\end{table}

\begin{table}[h]\center
\caption{Generators of the FJRW ring $\mathscr{H}_{W^T}, \quad W\in E^{(1,1)}_7$}
\label{ring-E7}
\begin{tabular}{|c|c|c|cc|c|}
\hline
$W^T$ &$G_{W^T}$& $\Theta(\be_{i,j,k})$&$\rho_1$&$\rho_2$&$\rho_{-1}$  \\
\hline

$X_1^4+X_2^4+X_3^2$ &$\mu_{4}^2\times\mu_{2}$& $\be[\frac{i}{4}],\be[\frac{j}{4}],\be[\frac{1}{2}]$&$\be_{2,1}$&$\be_{1,2}$&$\be_{2,2}=\rho_1^2\rho_2^2$  \\

\hline

$X_1^3X_2+X_1X_2^3+X_3^2$ &$\mu_{8}\times\mu_{2}$& $\be[-\frac{3i}{8}],\be[\frac{i}{8}],\be[\frac{1}{2}]$&$\be_{1}$&$\be_{5}$&$\be_{6}=\rho_1^2\rho_2^2$  \\
\hline

$X_1^3+X_1X_2^4+X_3^2$ &$\mu_{12}\times\mu_{2}$& $\be[\frac{-i}{3}],\be[\frac{i}{12}],\be[\frac{1}{2}]$&$\be_{1}$&$\be_{5}$&$\be_{10}=\rho_1\rho_2^3$  \\

\hline

$X_1^3+X_1X_2^2+X_2X_3^2$ &$\mu_{12}$& $\be[\frac{i}{3}],\be[-\frac{i}{6}],\be[\frac{i}{12}]$&$\be_{5}$&$\be_{1}$&$\be_{8}=\rho_1\rho_2\rho_3$  \\

\hline

\end{tabular}
\end{table}

\begin{table}[h]\center
\caption{Generators of the FJRW ring $\mathscr{H}_{W^T}, \quad W\in E^{(1,1)}_8$}
\label{ring-E8}
\begin{tabular}{|c|c|c|cc|c|}
\hline
$W^T$ &$G_{W^T}$& $\Theta(\be_{i,j,k})$&$\rho_1$&$\rho_2$&$\rho_{-1}$  \\
\hline

$X_1^6+X_2^3+X_3^2$ &$\mu_{6}\times\mu_{3}\times\mu_2$& $\be[\frac{i}{6}],\be[\frac{j}{3}],\be[\frac{1}{2}]$&$\be_{2,1}$&$\be_{1,2}$&$\be_{5,2}=\rho_1^4\rho_2$  \\

\hline

$X_1^3+X_1X_2^2+X_3^3$ &$\mu_{6}\times\mu_{3}$& $\be[\frac{-i}{3}],\be[\frac{i}{6}],\be[\frac{j}{3}]$&$\be_{1,1}$&$\be_{2,2}$&$\be_{4,2}=\rho_1\rho_2\rho_3$  \\
\hline

$X_1^4+X_1X_2^3+X_3^2$ &$\mu_{12}\times\mu_{2}$& $\be[\frac{-i}{4}],\be[\frac{i}{12}],\be[\frac{1}{2}]$&$\be_{2}$&$\be_{7}$&$\be_{9}=\rho_1^2\rho_2^2$  \\

\hline

\end{tabular}
\end{table}

\subsubsection{Classification of ring structure}\label{sec-ring}
According to isomorphism \eqref{ring-iso}, in order to classify the FJRW rings for ISES, it is enough to classify the Jacobian algebras.
For any special limit in the Saito-Givental theory for ISES, we need the ring structure for the classification.  Let us first focus on special limits at $\si=0$.
According to Saito \cite{S2}, simple elliptic singularities are classified
by their Milnor number and the elliptic curve at infinity. It follows
that the Jacobian algebras of the ISES with 3 variables can be classified
into $6$ isomorphic classes, parametrized by the pair consisting of the Milnor number
$\mu=\dim\mathscr{Q}_W$ and $j(0)$, the {\em j}-invariant of $E_{\si=0}$:

\begin{table}[h]\center
\caption{Classification of ring structure}
\label{table-ring}
\begin{tabular}{|l|l|}
\hline
$\Big(\mu,j(0)\Big)$&$W$\\
\hline
(8,0)&$X_1^3+X_2^3+X_3^3, X_1^2X_2+X_2^3+X_3^3, X_1^2X_2+X_1X_2^2+X_3^3, X_1^2X_2+X_2^2X_3+X_1X_3^2$\\
\hline
(8,1728)&$X_1^2X_2+X_2^2X_3+X_3^3$\\
\hline
(9,0)&$X_1^3X_2+X_2^4+X_3^2, X_1^3X_2+X_2^2X_3+X_3^2$\\
\hline
(9,1728)&$X_1^4+X_2^4+X_3^2, X_1^2X_2+X_2^2+X_3^4,X_1^3X_2+X_1X_2^3+X_3^2$\\
\hline
(10,0)&$X_1^6+X_2^3+X_3^2,X_1^3X_2+X_2^2+X_3^3$\\
\hline
(10,1728)&$X_1^4X_2+X_2^3+X_3^2$\\
\hline
\end{tabular}
\end{table}
For any two polynomials in the same list, it is easy to find a linear
map between the generators $X_1,X_2,X_3$ of the corresponding Jacobian algebras, such that
it induces a ring isomorphism. Let us point out that the choice of
such linear maps is not unique in general. 
In Section \ref{proof-t1}, we can always
adjust some constants such that the ring isomorphism will be extended
to an isomorphism of Frobeniu manifold, as well as an isomorphism of
the corresponding ancestor total potential. 

\subsubsection{Reconstruction of all genera FJRW invariants}\label{sec:FJRW-recons}
For an ISES $W^T$, its total ancestor potential $\mathscr{A}^{\rm FJRW}_{W^T}$ can be reconstructed from genus-0 primary correlators.   
This technique is already used in \cite{KS} for three special examples of ISESs. As the reconstruction procedures used there only require tautological relations on cohomology of moduli spaces of curves and the FJRW ring structure, we can easily generalize to all other examples. We sketch the general procedures here and refer to \cite{KS} for readers who are interested in more details. There are three steps.

First, we express the correlators of genus at least $2$ and the
correlators with descendant insertions in terms of correlators of genus-$0$ or genus-$1$ with
non-descendant insertions (called {\em primary correlators}). This step is based
on a tautological relation which splits a polynomial of $\psi$-classes
and $\kappa$-classes with higher degree to a linear combination of
products of boundary classes and polynomials of $\psi$-classes and
$\kappa$-classes of lower degrees. This is called
$g$-reduction. 
The reason why $g$-reduction works in our case is that the Selection
rule imposes a constraint on the degree of the polynomials involving
$\psi$- and $\kappa$- classes (see Theorem 6.2.1 in \cite{FJR}).
In general, for an arbitrary CohFT this
argument fails and one has to use other methods (e.g. Teleman's
reconstruction theorem). 

Next, we reconstruct the non-vanishing genus-$1$  primary correlators
from genus $0$ primary correlators using Getzler's relation. The latter
 is a relation in $H^4(\overline{\M}_{1,4})$, which gives identities involving
 the FJRW corrletors with genus $0$ and $1$. In order to obtain the
 desired reconstruction identity, i.e., to express genus-1 in terms of
 genus-0 correlators, one has to make an appropriate choice of the insertions corresponding to
 the 4 marked points in $\overline{\M}_{1,4}$ (see Theorem 3.9 in \cite{KS}).

Finally, to reconstruct the genus-$0$ correlators we use the WDVV equations. 
We say that a homogeneous element $\alpha\in \mathscr{H}_{W^T}$ is \emph{primitive} if it
cannot be decomposed as a product $a'\bullet a''$ of two
elements $a'$ and $a''$ of non-zero degrees.
We also say that a genus-0 correlator is a \emph{basic correlator} if there are at most
two non-primitive insertions, neither of which is the identity. We use
the WDVV equation to rewrite a primary genus-$0$ correlator which contains
several non-primitive insertions to correlators with fewer
non-primitive insertions and correlators with a fewer number of marked
points. Again the Selection rule should be taken into account in order
to obtain a bound for the number of
marked points. It turns out that all correlators are
determined by the basic correlators with at most four marked
points (see Lemma 3.7 in \cite{KS}).
\begin{lemma}\label{FJRW-recons}
For an invertible simple elliptic singularity $W^T$ the total ancestor
FJRW potential $\mathscr{A}_{W^T}^{\rm FJRW}$ of $(W^T,G_{W^T})$ is reconstructed from the pairing, the FJRW
ring structure constants and the $4$-point basic correlators with one of
the insertions being a top degree element.
\end{lemma}

\subsubsection{The 4-point genus-0 FJRW invariants}
We introduce the following notation.
\begin{definition}\label{xi}
Let $\Xi(\rho_1,\rho_2,\rho_3)$ be a degree 1 monomial with leading
coefficient 1. For simplicity, we denote by $\langle \Xi,
\rho_{-1}\rangle_{0,4}^{W^T}$ a basic correlator such that the first
three insertions give a factorization of $\Xi$. 
\end{definition}
This is well defined because the WDVV equations guarantee that $\langle \Xi,\rho_{-1}\rangle_{0,4}^{W^T}$ does not depend on the choices of the factorization. For example, let
$\Xi(\rho_1,\rho_2,\rho_3)=\rho_1^2\rho_2^2$; then the notation
$\langle\Xi, \rho_{-1}\rangle_{0,4}^{W^T}$ represents any of the
following choices of correlators:
$$\LD\rho_1,\rho_1,\rho_2^2, \rho_{-1}\RD_{0,4}^{W^T},\LD\rho_1,\rho_2,\rho_1\rho_2, \rho_{-1}\RD_{0,4}^{W^T},\LD\rho_2,\rho_2,\rho_1^2, \rho_{-1}\RD_{0,4}^{W^T}.$$
But it does not represent $\LD 1,\rho_1,\rho_1\rho_2^2, \rho_{-1}\RD_{0,4}^{W^T}$, which is not a basic correlator.
\begin{lemma}\label{lm:recons-FJRW}
Let $W^T$ be an ISES; then the total FJRW potential $\mathscr{A}_{W^T}^{\rm FJRW}$ for $(W^T, G_{W^T})$
can be reconstructed from the FJRW algebra, and the basic 4-point FJRW
correlators $\langle \Xi,\rho_{-1}\rangle_{0,4}^{W^T}$.
Furthermore, if $W^T$ is an ISES as in Tables \ref{ring-E6},
\ref{ring-E7}, or \ref{ring-E8}, then
\begin{equation}\label{eq:4pt}
\LD\Xi(\rho_1,\rho_2,\rho_3), \rho_{-1}\RD_{0,4}^{W^T}
=
\begin{cases}
q_i^T & \mbox{ if } \Xi=M_i,\\
0 & \mbox{ otherwise. }
\end{cases}
\end{equation}
where $M_i$ are the homogeneous monomials such that $W=M_1+M_2+M_3$. 
\end{lemma}
\proof
In \cite{KS} Theorem 3.4, it was proved for three special simple elliptic singularities $W^T=X_1^3+X_2^3+X_3^3, X_1^3+X_1X_2^2+X_2X_3^2$ and $X_1^3+X_1X_2^2+X_3^3$, their FJRW correlators with symmetry group $G_{W^T}$ can be reconstructed from their FJRW algebra and some basic 4-point correlators. We apply the same method to all cases of simple elliptic singularities here. Finally, using WDVV equations in each case, it is again not hard to verify all 4-point basic correlators without insertion $\rho_{-1}$ can be reconstructed too.

For the second part of the lemma, we use WDVV and concavity to compute FJRW correlators.
We show the argument works for singularities of Fermat type and of loop type. Other cases are similar. 
For a Fermat type singularity, put $M_i=X_i^{1/q_i^T}$, 
since all insertions are
narrow, we apply the Concavity Axiom \eqref{eq:OGRR} to compute
\beq\label{cor}
\LD
\rho_i,\rho_i,\rho_i^{1/q_i^T-2}, \rho_{-1}
\RD_{0,4}^{W^T}.
\eeq  
Note that $\deg\mathscr{L}_i=-2$ and the degree shifting numbers
are $(2q_i^T,2q_i^T,1-q_i^T,1-q_i^T)$, thus the dual graphs will have
$\Theta_{\Gamma}
=0,0,1-3q_i^T$. The correlator \eqref{cor} becomes
\ben
\frac{1}{2}\left(B_2(q_i^T)+B_2(1-3q_i^T)+2B_2(0)-2B_2(q_i^T)-2B_2(1-q_i^T)\right)=
q_i^T.
\een
For loop type, $W^T=X_1^2X_3+X_1X_2^2+X_2X_3^2$. Let us compute $\LD\rho_1,\rho_1,\rho_2,\rho_{-1}\RD_{0,4}^{W^T}$, which is not concave. However, the Concavity Axiom \eqref{eq:OGRR} implies 
$$\LD\be_2,\be_4,\be_7,\be_2\RD_{0,4}^{W^T}=-\frac{2}{9}.$$
On the other hand, WDVV equations show
$$\left\{
\begin{array}{l}
\LD\be_1\bullet\be_4,\be_4,\be_7,\be_2\RD_{0,4}^{W^T}+\LD\be_1,\be_4,\be_4,\be_7\bullet\be_2\RD_{0,4}^{W^T}=\LD\be_7\bullet\be_4,\be_4,\be_1,\be_2\RD_{0,4}^{W^T};\\
\LD\be_4\bullet\be_4,\be_1,\be_7,\be_2\RD_{0,4}^{W^T}+\LD\be_4,\be_4,\be_1,\be_7\bullet\be_2\RD_{0,4}^{W^T}=\LD\be_4\bullet\be_7,\be_1,\be_4,\be_2\RD_{0,4}^{W^T};\\
\end{array}
\right.$$
We observe up to symmetry, $\LD\be_5,\be_1,\be_7,\be_2\RD_{0,4}^{W^T}=\LD\be_8,\be_1,\be_4,\be_2\RD_{0,4}^{W^T}$. Recall the ring relations in Example \ref{ex:loop}, we obtain
\ben
\LD\rho_1,\rho_1,\rho_2,\rho_{-1}\RD_{0,4}^{W^T}=\LD\be_4,\be_4,\be_1,\be_6\RD_{0,4}^{W^T}=\frac{1}{3}.
\qed
\een

\section{\texorpdfstring{Mirror symmetry at $\si=0$}{Zero}}\label{MS-0}

In this section, our goal is to prove Theorem \ref{t1}. According to the reconstruction
results in FJRW theory (see Lemma \ref{lm:recons-FJRW}) we need to
compute certain 3- and 4- point genus-0 correlators in Saito's theory
and compare them to the ones in the mirror FJRW
theory. 

\subsection{The Saito--Givental limit}\label{SG-CohFT}
The higher-genus Saito--Givental CohFT is in general defined only at
semisimple points $\s\in S$, because the asymptotic operator $R_\s(z)$ (see
Section \ref{sec:givental}) has singularities along the {\em
  caustic} $\mathcal{K}\subset \S$ consisting of non-semisimple $\s$.  However, the
genus-0 part of the Saito--Givental CohFT is well defined for all
$\s\in \S$. The logic of our argument is the following: we identify first the genus-0
CohFTs; then we use that the higher-genus reconstruction of the
correlation functions (see Section \ref{sec:FJRW-recons}) works for both theories. There are two conclusions from
this. First, the FJRW total ancestor potential is convergent at all semisimple $\s\in
\S$, such that $\si=s_{\mu-1}$ is sufficiently close to 0 and it coincides
with the Saito--Givental ancestor potential $\mathscr{A}^{\rm SG}_W(\s)$ (see \eqref{SG-ancestor}). The second conclusion is that 
$\mathscr{A}^{\rm SG}_W(\s)$ extends holomorphically through the caustic $\mathcal{K}$. 
In other words, one can define the Saito-Givental limit $\mathscr{A}_W^{\rm SG}(\sigma)$ near point $\si=0$ by
\beq\label{limit-SG}
\mathscr{A}_W^{\rm SG}(\sigma):=\lim_{\s\to(\si, {\bf 0})}\mathscr{A}^{\rm SG}_W(\s).
\eeq

There is an alternative way to proceed provided that we know that the
genus-0 CohFTs are the same. It is based on Teleman's
classification of semisimple CohFTs \cite{Te}. We refer to \cite{CI} and Section 5
in \cite{MRS} for more details. The advantage of this approach is that
one can prove a stronger result. Namely, the higher-genus
Saito--Givental CohFT (not only its ancestor potential) extends
through the caustic $\mathcal{K}$. 

The details in the proofs of the above statements can be found in
Lemma 3.2 in \cite{MR} for the case of ancestor potentials only and in 
Proposition 5.5 in \cite{MRS} for the CohFTs. Recently, the extension problem is solved for generic isolated singularityby the first author \cite{M}, using Eynard-Orantin recursion.

\subsection{B-model 3-point genus-0 correlators}
We continue to use the same notation as in Section \ref{PF:system}. Namely,
let $W=M_1+M_2+M_3$ be an ISES with a miniversal deformation given by a monomial
$\phi_\m(\x)$, $\m=(m_1,m_2,m_3)$. 
We choose a primitive form $\omega=d^3\x/\pi(\si)$ in a neighborhood
of $\si=0$, such that $\pi(\si)$
is the solution to the Picard-Fuchs equation \eqref{PF:1} satisfying
the initial conditions $\pi(0)=c,\pi'(0)=0$, where the constant $c$ is
such that the residue pairing (see \eqref{res:pairing}) satisfies
\ben
\left. \langle 1,\phi_\m\rangle\right|_{\s=0} =1\, .
\een
Let $\{t_{\r}\}$ be the flat
coordinate system, such that $t_\r(0)=0$ and the flat vector fields $\d_\r:=\d/\d t_\r$
agree with $\d /\d s_\r$ at $\s=0$.  

The primitive form induces an isomorphism between the tangent and the
vanishing cohomology bundle via the following period mapping:
\beq\label{period-map}
\d/\d t_\r\mapsto -\nabla_{\frac{\d}{\d \lambda}}^{-1}\,
\nabla_{\frac{\d}{\d t_\r}}\, \int \frac{\omega}{d F} = \int \delta_\r(\s,\x)
\, \frac{\omega}{d F}=\int \delta_\r(\s,\x)
\, \frac{1}{\pi(\si)}\,\frac{d^3\x}{d F},
\eeq
where $\delta_\r$ is some homogeneous polynomial (in $\x$) of degree ${\rm
  deg}(\phi_\r)$. Note that the Kodaira--Spencer
isomorphism takes the form
\beq\label{iso:KS}
\d/\d t_\r\mapsto \delta_\r(\s,\x) \ {\rm mod } \ (F_{X_1},F_{X_2},F_{X_3}).
\eeq
We know $\d_{\r_1}\bullet\d_{\r_2}$ is induced from multiplication in $q_*\mathcal{O}_C$, and the pairing is
$$\LD\delta_{\r_1},\delta_{\r_2}\RD:=\eta(\d_{\r_1},\d_{\r_2})={\rm Res}\frac{\delta_{\r_1}(\s,\x) \delta_{\r_2}(\s,\x)}{(\d_{X_1}W_\si)\, (
  \d_{X_2}\,W_\si) (\d_{X_3}W_\si)}\frac{d^3\x}{\pi(\si)^2}.$$
By definition, the restriction of the 3-point correlators to the
marginal direction is
\beq\label{3pt}
\LD \delta_{\r_1},\delta_{\r_2},\delta_{\r_3}\RD_{0,3} 
=\LD \delta_{\r_1},\delta_{\r_2}\cdot\delta_{\r_3}\RD
= {\rm Res}\,
\frac{\delta_{r_1}(\si,\x)\delta_{r_2}(\si,\x)\delta_{r_3}(\si,\x)}{(\d_{X_1}W_\si)\, (
  \d_{X_2}\,W_\si) (\d_{X_3}W_\si) }\, \frac{d^3\x}{\pi(\si)^2}.
\eeq
Note that the 3-point correlator depends only on the product $\Xi:=
\delta_{\r_1}\delta_{\r_2}\delta_{\r_3}.$ Therefore we can simply use
the notation $\langle\Xi\rangle_{0,3}$ instead. Finally, definition
\eqref{3pt} makes sense even if we replace $\delta_\r$,
$\r=\r_1,\r_2,\r_3$ by arbitrary polynomials, not only the ones that
correspond to flat vector fields via \eqref{iso:KS}.

\subsection{B-model 4-point genus-0 correlators}
Let $\mathscr{F}_{0}^{\rm SG}$ be the genus-0 generating functions for the Frobenius manifold of miniversal deformations near the origin. By definition,
\beqa
\LD\delta_{\r_1},\dots,\delta_{\r_n}\RD_{0,n}^{\rm SG}=\frac{\partial^n\mathscr{F}_{0}^{\rm SG}}{\partial\,t_{\r_1}\dots\partial\,t_{\r_n}}\Big\vert_{\t=0}.
\eeqa
Thus, using that $\d/\d\si=\delta_\m$ at $\si=0$, we get
\beq\label{SG4pt}
\LD\delta_{\r_1},\delta_{\r_2},\delta_{\r_3},\delta_{\m}\RD_{0,4}^{\rm SG}
=\d_\si\,\LD\delta_{\r_1},\delta_{\r_2},\delta_{\r_3}\RD\Big\vert_{\si=0}\ .
\eeq
In order to compute 4-point correlators of the form \eqref{SG4pt} it
is enough to determine $\delta_\r(\si,\x)$ up to linear terms in $\si$. 
To begin with, we notice that $\phi_{\r+\m}$ lies in the Jacobian
ideal of $W_{\sigma}$. More precisely, the following Lemma holds.
\begin{lemma}\label{Jacobi-rel}
There are polynomials $g_{\r,i}\in \C[\si,X_1,X_2,X_3]$ such that 
$$(1-C\si^l)\ \phi_{\r+\m}=\sum_{i=1}^{3}g_{\r,i}\,\d_iW_{\sigma}.$$
\end{lemma}
This Lemma can be proved in all cases by using Saito's higher residue
pairing. However, in what follows, we need an explicit formula for 
$$g_{\r}:=\left(g_{\r,1},g_{\r,2},g_{\r,3}\right).$$ 
Therefore we verified the lemma on a case-by-case basis. Some of our computations will be
given below. The remaining cases are completely analogous.

There are several corollaries of Lemma \ref{Jacobi-rel}. First of
all, note that under the period map \eqref{period-map} the
Gauss--Manin connection takes the form \eqref{frob_eq1} (with $z\equiv
-\d_\lambda^{-1}$. It follows that if ${\rm deg}(\phi_\r)$ is not
integral, then the restriction of the section \eqref{period-map} of
the vanishing cohomology bundle to the marginal deformation subspace
must be flat, i.e., the sections
\beq\label{s-geom}
[\delta_\r\,\omega](\si):=\int\,\delta_\r(\si,\x)\,\frac{\omega}{d W_{\si}},\quad {\rm deg}(\phi_\r)\notin\Z
\eeq
are independent of $\si$. 
Furthermore, using formulas \eqref{GM-con} for the Gauss-Manin
connection we get
$$
(1-C\si^l)\frac{\d}{\d\sigma}\Phi_{\r}=-\int \sum_{i=1}^{3}\d_i\,
g_{\r,i}\, \frac{d^3\x}{d W_\si}
$$
Both sides must have the same degree, i.e., 
\beq\label{aux-eq}
(1-C\si^l)\frac{\d}{\d\sigma}\Phi_{\r}=\sum_{\r'}c_{\r,\r'}(\sigma)\Phi_{\r'}, 
\eeq
where the sum is over all $\r'$, such that $\deg\phi_{\r}=\deg\phi_{\r'}$
and $c_{\r,\r'}(\sigma)\in \C[\si]$ are some polynomials. 
\begin{lemma}\label{lin-approx}
Suppose ${\rm deg}(\phi_\r)\notin \Z$; then we have 
\begin{equation}\label{eq:first}
\delta_{\r} 
=\phi_{\r}-\sigma\sum_{\r',\r'\neq\r}c_{\r,\r'}(0)\,\phi_{\r'} + O(\si^2),
\end{equation}
where $O(\si^2)$ denotes terms that have order of vanishing at $\si=0$ at least 2.
\end{lemma}
\proof
Follows easily from \eqref{aux-eq}. We omit the details.
\qed

Let $M(X_1,X_2,X_3)\in \C[\x]$ be a weight-1 monomial with leading coefficient
1. Our next goal is to evaluate the following auxiliar expression (recall Definition \ref{xi}):
\ben
\langle M,\phi_\m\rangle_{0,4}:=\left. \d_\si\langle M\rangle_{0,3}\right|_{\si=0}.
\een
\begin{lemma}\label{lm:4pt=0}
The number $\langle M,\phi_\m\rangle_{0,4}$ is non-zero iff $M=M_i$ for
some $i=1,2,3$. In the latter cases the numbers are given as follows
\beq\label{SG3pt}
\Big(\langle M_1,\phi_\m\rangle_{0,4},\langle M_2,\phi_\m\rangle_{0,4},\langle M_3,\phi_\m\rangle_{0,4}\Big)
=-(m_1,m_2,m_3)E_{W}^{-1}.
\eeq
\end{lemma}
\proof
For the second part, we apply the operators $X_i\, \d_{X_i}$, $i=1,2,3,$
to the identity
\ben
M_1+M_2+M_3 = W_\si-\si\phi_\m(\x)
\een
and take the residue. We get 
\ben
\langle M_1\rangle_{0,3}\, a_{1i} + \langle M_2\rangle_{0,3}\, a_{2i}
+ \langle M_3\rangle_{0,3}\, a_{3i} =
-\si\,m_i\, \langle \phi_\m\rangle_{0,3}.
\een
It remains only to differentiate with respect to $\si$ and set
$\si=0$.

For the first part, because $M$ is a weight-$1$ monomial with coefficient $1$, we can use the
relations in the Jacobian algebra of $W_{\si}$ to rewrite $M$ as a
product of  $\phi_\m$ and a function of $\sigma$. Let us write $M=h(\sigma)\,\phi_{m}$. 
For example, in the Fermat $E_{6}^{(1,1)}$ case, 
$$X_1^3=-3\si\phi_{111};\quad (1+\frac{\si^3}{27})\,X_1^2X_2=0.$$ 
If $M\neq M_i$, $i=1,2,3$, then $h(\si)$ either does not vanish at $\si=0$ or vanishes at $\si=0$ with order at least $2$. In both cases,  $\langle M,\phi_\m\rangle_{0,4}$ vanish.
\qed

Now we are ready to compute the 4-point correlators that are needed
for the reconstruction of the CohFT. Let $\delta_{\r}(\s,x)$,
$\r=\r_1,\r_2,\r_3$ be polynomials corresponding to the flat vector
fields $\d/\d t_\r$ via the Kodaira--Spencer isomorphism
\eqref{iso:KS}. Put
\ben
\Xi(\s,\x)=\delta_{\r_1}(\s,x) \delta_{\r_2}(\s,x) \delta_{\r_3}(\s,x).
\een
Note that $\Xi(0,\x)$ is a homogeneous monomial (see \eqref{eq:first}) with leading coefficient 1. 
\begin{lemma}\label{lm:SG4pt}
The $4$-point genus-$0$ correlators with a top degree insertion $\delta_\m$ are 
\ben
\LD \delta_{\r_1},\delta_{\r_2},\delta_{\r_3} , \delta_\m\RD_{0,4}^{\rm SG}
=
\begin{cases}
-q_i^T &\mbox{, if } \Xi(0,\x)=M_i,\\
\ \ 0        & \mbox{, otherwise.}
\end{cases}
\een
\end{lemma}
\proof
The same argument as in the proof of the first part of Lemma \ref{lm:4pt=0} also works for $\Xi(\sigma,\x)$.  
Thus if $\Xi\neq M_i$,
i=1,2,3, we have 
\ben
\LD\Xi,\delta_\m\RD_{0,4}^{\rm SG}=0.
\een
In order to finish the proof we need only to compute the correlators
when $\Xi(0,\x)=M_i$ for some $i=1,2,3$. Note that the diagonal
entries of the matrix $E_W$ are always at least 2 (see Table
\ref{table-1}). Therefore, it is enough to compute the following correlators:
$$
\left\{
\begin{array}{l}
\LD\delta_{100},\delta_{100},\delta_{\r},\delta_\m\RD_{0,4}^{\rm SG},\quad
\r=(a_{11}-2,a_{12},a_{13}),\\
\LD\delta_{010},\delta_{010},\delta_{\r},\delta_\m\RD_{0,4}^{\rm SG},\quad
\r=(a_{21},a_{22}-2,a_{23}),\\
\LD\delta_{001},\delta_{001},\delta_{\r},\delta_\m\RD_{0,4}^{\rm SG}, \quad 
\r=(a_{31},a_{32},a_{33}-2).
 \end{array}
\right.
$$
We do not have a uniform computation since we need to use Lemma
\ref{lin-approx}, for which the coefficients $c_{\r,\r'}(0)$ can be
computed only on a case-by-case basis. Let us sketch the main steps of
the computation in several examples, leaving the details and the
remaining cases to the reader. We will make use of the notation 
\ben
\delta(\si,\x)\approx \phi(\si,\x),\quad \delta,\phi\in \C[\x],
\een
which means first order approximation at $\si=0$, i.e.,
$\delta(\si,\x)-\phi(\si,\x)=O(\si^2)$.

{\em Case 1:}
$W=X_1^3+X_2^3+X_3^3\in E_{6}^{(1,1)}$ and $\phi_\m=X_1X_2X_3$. Since $W$ is symmetric in
$X_1,X_2,X_3$ it is enough to compute only one of the correlators, say
$\Xi=M_1$. After a straightforward computation (the notation is the
same as in Lemma \ref{Jacobi-rel}) we get
\ben
g_{100}=\Big(\frac{1}{3}\, \phi_{011},-\frac{\si}{9}\,
\phi_{002},\frac{\si^2}{27}\, \phi_{101}\Big).
\een
It follows that $\delta_{100}\approx\phi_{100}$ and then using formula
\eqref{SG3pt} we get
\ben
\LD\delta_{100},\delta_{100},\delta_{100},\delta_\m\RD_{0,4}^{\rm SG}=-\frac{1}{3}.
\een
{\em Case 2:}
$W=X_1^4+X_2^4+X_3^2\in E_{7}^{(1,1)}$ and $\phi_\m=X_1^2X_2^2$.
In this case $M_3=0$ in the Jacobian algebra of $W$ and $W$ is symmetric in
$X_1$ and $X_2$. It is enough to compute only one of the correlators,
say the one with $\Xi(0,\x)= M_1.$ We have
\ben
g_{100}=\Big(\frac{1}{4}\, \phi_{020},-\frac{\si}{8}\, \phi_{110},0\Big).
\een
It follows that 
\ben
\delta_{100}\approx\phi_{100},\quad\delta_{200}\approx\phi_{200}+\frac{\si}{4}\phi_{020}.
\een
Using formula \eqref{SG3pt} we find
$$
\LD\delta_{100},\delta_{100},\delta_{400},\delta_\m\RD_{0,4}^{\rm SG}=-\frac{1}{4}.
$$
{\em Case 3:}
$W=X_1^3X_2+X_2^2+X_3^3\in E_{8}^{(1,1)}$ and $\phi_\m=X_1X_2X_3.$ In
this case, since $M_2=0$ in the Jacobian algebra, we need to compute two
correlators. We have
$$
\left(\begin{array}{c}
g_{100}\\
g_{010}\\
g_{001}
\end{array}\right)
=
\left(
\begin{array}{ccc}
\frac{1}{3}\phi_{001}-\frac{\si^2}{54}\phi_{200}&\frac{\si^2}{18}\phi_{110}&-\frac{\si}{9}\phi_{010}\\
-\frac{1}{6}\phi_{201}+\frac{\si^2}{27}\phi_{110}&\frac{1}{2}\phi_{111}&-\frac{\si^2}{9}\phi_{210}\\
-\frac{\si}{9}\phi_{010}-\frac{\si^2}{54}\phi_{101}&\frac{\si^2}{9}\phi_{011}&\frac{1}{3}\phi_{110}
\end{array}
\right).
$$
It follows that we have the following linear approximations:
$$\delta_{100}\approx\Phi_{100},\quad\delta_{001}\approx\Phi_{001},\quad\delta_{110}\approx\Phi_{110}.$$
The correlators then become
$$\LD\delta_{100},\delta_{100},\delta_{110},\delta_\m\RD_{0,4}^{\rm SG}=-\frac{1}{3},\quad \LD\delta_{001},\delta_{001},\delta_{001},\delta_\m\RD_{0,4}^{\rm SG}=-\frac{1}{3}$$
{\em Case 4:}
The Fermat type $E_{8}^{(1,1)}$, i.e.  $W=X_1^6+X_2^3+X_3^2$ and
$\phi_\m=X_1^4X_2.$ In this case $M_3=0$, so again we have to compute two
correlators. We have
$$g_{100}=\Big(\frac{1}{6}\phi_{020}+\frac{\si^2}{27}\phi_{200},-\frac{2\si}{9}\phi_{300},0\Big), \quad 
g_{010}=\Big(-\frac{\si}{18}\phi_{300}+\frac{\si^2}{27}\phi_{110},\frac{1}{3}\phi_{400},0\Big).$$
It follows that the first order approximations that we need are
$$
\delta_{100}\approx\Phi_{100},\quad \delta_{010}\approx\Phi_{010},\quad
\delta_{400}\approx\Phi_{400}+\frac{\si}{2}\Phi_{210}.
$$
Formulas  \eqref{SG3pt} and \eqref{SG4pt} imply
$$
\LD\delta_{100},\delta_{100},\delta_{400},\delta_\m\RD_{0,4}^{\rm SG}=-\frac{1}{6};\quad
\LD\delta_{010},\delta_{010},\delta_{010},\delta_\m\RD_{0,4}^{\rm SG}=-\frac{1}{3}.
\qed
$$

\subsection{Proof of Theorem \ref{t1}}\label{proof-t1}
The 4-point correlators in Lemma \ref{lm:recons-FJRW} and
Lemma \ref{lm:SG4pt} have opposite signs. If $\rho_{-1}\mapsto\delta_{\m}$, we rescale the primitive form by $(-1)$ and define 
\beq\label{mirror-Gepner}
\mathscr{H}_{W^T}\to T^*\S_W,\quad 
\rho_{\r} \mapsto (-1)^{1-\deg\phi_{\r}}\delta_{\r},\quad
\r=(r_1,r_2,r_3),
\eeq
where $\rho_\r=\rho_1^{r_1}\rho_2^{r_2}\rho_3^{r_3}$. For the new basis, the
3-point correlators in the Saito--Givental theory do not change, while
the 4-point correlators are rescaled by $(-1)$.  
Lemma \ref{lm:recons-FJRW} and Lemma
\ref{lm:SG4pt} imply that the map \eqref{mirror-Gepner} identifies $\mathscr{A}^{\rm FJRW}_{W^T}$ and $\mathscr{A}^{\rm SG}_{W}(\sigma=0)$. If $\rho_{-1}\mapsto c\delta_{\m}$ for some nonzero constant $c>1$, we need to rescale the ring generators furthermore to obtain an suitable mirror map.
Thus Theorem \ref{t1},
a) is proved.


On the other hand, since we already computed all basic 4-point genus-0 correlators for Saito-Givental limit at $\si=0$, we can use it here to identify two special limit points of Saito-Givental theory, if they have isomorphic rings. We notice that all the ring structures at $\si=0$ are already listed in Section \ref{sec-ring}. The following lemma gives a complete classification for $\si=0$, in the sense of Definition \ref{def-special-limit}.  In particular it gives a proof for part b) of Theorem \ref{t1}.
\begin{lemma}\label{FJRW-rings}
If $W_1$ and $W_2$ are invetible simple elliptic singularities, $\mathscr{Q}_{W_1}\cong\mathscr{Q}_{W_2}$. Then there exists a ring isomorphicm $\Psi:\mathscr{Q}_{W_1}\cong\mathscr{Q}_{W_2}$, such that
$\Psi$ preserves the Saito-Givental limits at $\si=0$, i.e. $\Psi:\mathscr{A}^{\rm SG}_{W_1}(\si)=\mathscr{A}^{\rm SG}_{W_2}(\si)$ for $\si=0$.
\end{lemma}
\proof
We will construct explicitly linear
isomorphisms $\Psi$ inducing the ring
isomorphisms; then one has to check that they also preserve the $4$-point correlators in Lemma \ref{lm:SG4pt}. We list construction of isomorphisms $\Psi$ for some examples below. We will explain the notations through one example later.
For each $W$, we can fix a particular choice of $\phi_\m$ in the table and choose $\{Y_i=X_i/\la_i\}_{i=1}^3$. For all other choices of $\phi_\m$, we can obtain an isomorphism $\Psi$ by rescaling the ring generators appropriately.
\begin{table}[h]\center
\caption{Classification of $\si=0$}
\center
\begin{tabular}{|c|c|c|c|}
\hline
W$$&$\phi_\m$& {\rm Ring \ Generators} & {\rm Constraints}\\
\hline
\hline
$X_1^3+X_2^3+X_3^3$ &$X_1X_2X_3$& $(Y_1, Y_2,Y_3)$ & \\
\hline
$X_1^2X_2+X_2^3+X_3^3$& $X_2^2X_3$ & $(\frac{X_1}{\sqrt{3}}+X_2,\frac{-X_1}{\sqrt{3}}+X_2,X_3)$ & $\la_1^4\la_2\la_3=\la_1\la_2^4\la_3=\frac{\la_1\la_2\la_3^3}{2}=\frac{1}{4}$\\
\hline
$X_1^2X_2+X_1X_2^2+X_3^3$&$X_1X_2X_3$&$(\frac{X_1}{\be[\frac{5}{6}]}+X_2,\frac{X_1}{\be[\frac{1}{6}]}+X_2,X_3)$&$\la_1^4\la_2\la_3=-\la_1\la_2^4\la_3=\frac{\la_1\la_2\la_3^4}{\sqrt{-27}}=\frac{1}{\sqrt{-27}}$\\
\hline
\hline
$X_1^3X_2+X_2^4+X_3^2$ &$X_1X_2^3$& $(Y_1,Y_2)$  & \\
\hline
$X_1^3X_2+X_2^2X_3+X_3^2$&$X_1X_2X_3$&$(X_1,X_2)$&$-2\la_1^4\la_2^4=8\la_1\la_2^7=1$\\
\hline
\hline
$X_1^4+X_2^4+X_3^2$ &$X_1^2X_2^2$& $(Y_1,Y_2)$ & \\
\hline
$X_1^2X_2+X_2^2+X_3^4$&$X_1^2X_3^2$&$(X_1,X_3)$&$-4\la_1^6\la_2^2=\la_1^2\la_2^6=1$\\
\hline
$X_1^3X_2+X_1X_2^3+X_3^2$&$X_1^2X_2^2$&$(X_1+X_2,-X_1+X_2)$&$-64\la_1^6\la_2^2=64\la_1^2\la_2^6=1$\\
\hline
\hline
$X_1^6+X_2^3+X_3^2$ &$X_1^4X_2$& $(Y_1,Y_2)$ &\\
\hline
$X_1^3X_2+X_2^2+X_3^3$&$X_1X_2X_3$&$(X_1,X_3)$&$8\la_1^{10}\la_2=-2\la_1^4\la_2^4=1$\\
\hline
\end{tabular}
\end{table}

Now we consider the example of $W_1=X_1^3+X_2^3+X_3^3$ and $W_2=X_1^2X_2+X_2^2X_3+X_1X_3^2$, both of $\phi_\m=X_1X_2X_3$. We construct a ring homomorphism
$\Psi:\mathscr{Q}_{W_1}\to\mathscr{Q}_{W_2}$
which sends a basis of generators of $\mathscr{Q}_{W_1}$, $\{Y_i=X_i/\la_i\}_{i=1}^3$ to $\{\Psi(Y_i)\}_{i=1}^{3}\in\mathscr{Q}_{W_2}$,
$$
\left(\Psi(Y_1),\Psi(Y_2),\Psi(Y_3)\right)
=
\left(
X_1,
X_2,
X_3
\right)
\left( \begin{array}{ccc}
1 &1 & 1 \\
1 & \be[\frac{1}{3}]& \be[\frac{2}{3}] \\
1 & \be[\frac{2}{3}] & \be[\frac{1}{3}] \end{array} \right).
$$
We choose parameters $\la_1,\la_2,\la_3$ satisfying
\beq\label{FJRW-iso}
\la_1^4\la_2\la_3=\be[\frac{1}{3}]\,\la_1\la_2^4\la_3=\be[\frac{2}{3}]\,\la_1\la_2\la_3^4=-\frac{1}{27}.
\eeq
We can check all the 4-point genus-0  correlators listed in Lemma \ref{lm:SG4pt} match under the  isomorphism $\Psi$ using \eqref{FJRW-iso}. For example, since both of
$\phi_{\m}=X_1X_2X_3$, we get 
$$\LD\Psi(X_1),\Psi(X_1),\Psi(X_1),\Psi(X_1X_2X_3)\RD_{0,4,W_2}
=-27\la_1^4\la_2\la_3\LD X_1,X_1,X_1,X_1X_2X_3\RD_{0,4,W_2}=-\frac{1}{3}.$$
Similarly, the other two 4-point genus-0 correlators match follows from the other two identities in \eqref{FJRW-iso}. Thus $\Psi:\mathscr{A}^{\rm SG}_{W_1}(\si)=\mathscr{A}^{\rm SG}_{W_2}(\si)$ for $\si=0$ according to Lemma \ref{lm:recons-FJRW} and Lemma
\ref{lm:SG4pt}.
\qed

\section{\texorpdfstring{Global mirror symmetry for Fermat simple elliptic singularities}{Unity}}\label{MS-punctures}
The goal in this section is to prove Theorem \ref{t2}. 
For special limit $\si\neq0, \infty$, it is enough to prove it only for one of the points
$p_k=C^{-1/l}\eta^k$, $1\leq k\leq l$. The
remaining cases follow easily due to the $\Z/l\,\Z$-symmetry of
$\Si$. For all those limits except $\sigma=\infty$ for $W=X_1^6+X_2^3+X_3^2$, according to the reconstruction Lemma \ref{thm:GW recons}, we
need to construct an appropriate {\em mirror map} from the Chen--Ruan orbifold cohomology ring to the limit of Jacobian algebras, such that after choosing an appropriate primitive form, the
Poincar\'e pairing is identified with the residue pairing and the
3-point correlator (see Lemma \ref{thm:GW recons} ) is the same for
both the Gromov--Witten and the Saito--Givental CohFTs. Finally, we also prove the Saito-Givental limit $\sigma=\infty$ for $W=X_1^6+X_2^3+X_3^2$ is isomorphic to an FJRW theory. This agrees with the physicists' prediction that the monodromy of the Gauss--Manin connection around
the large volume limit point should be maximally unipotent, while as we will
see below, the monodromy around $\si=\infty$ is diagonalizable.

The limit of the Saito-Givental
theory of ISESs at $\si=\infty$ is already discussed in
\cite{MR} for $(W,\phi_\m=X_1X_2X_3)$, where
\ben
W=X_1^3+X_2^3+X_3^3\in E_{6}^{(1,1)}, 
X_1^3X_2+X_2^2X_3+X_3^2\in E_{7}^{(1,1)},
X_1^3X_2+X_2^2+X_3^3\in E_{8}^{(1,1)}.
\een 
Namely, it was proved that the Saito-Givental theory at $\si=\infty$ is
mirror to the Gromov-Witten theory respectively of
$\mathbb{P}^1_{3,3,3}, \mathbb{P}^1_{4,4,2}$ and
$\mathbb{P}^1_{6,3,2}.$

\subsection{Construction of a mirror map}
We construct a mirror map based on solving the systems of hypergeometric equations. We will introduce explicit how to construct it near $\sigma=p_k$. For the special limit $\sigma=\infty$, the construction is similar.
 
\subsubsection{Non-twisted sectors}
The primitive form is chosen to be $\omega=d^3\x/\pi_A(\si)$, where
the cycle $A\in H_1(E_{\si})$ is invariant with respect to the local
monodromy around $\si=p_k$. Recall that $\pi_A$ is a solution to
the hypergeometric equation \eqref{HG:1} with weights $(\alpha,\beta,\gamma)$, where $\gamma=\alpha+\beta$. The invariance of $A$ implies near $x=1$, we have
\beq\label{lambda_W}
\pi_{A}(\si)=\la_W\,F_1^{(1)}(x):=\la_W\ _2F_1\left(\alpha,\beta,1,1-x\right), \quad \la_W\in\C^*.
\eeq
Since the $j$-invariant of $E_\si$ always has the form 
$$j(\si)=\frac{P(\si)}{(1-C\si^l)^N}$$ for some
polynomial $P(\si)\in\C[\sigma]$ and some integer $N$.
We can always choose a second cycle $B\in H_1(E_{\si})$, such that 
\beq\label{K-constant}
\pi_B(\si)=\frac{N\la_W}{2\pi\sqrt{-1}}\left(F_2^{(1)}(x)-\frac{\ln\, P(p_k)}{N}\, F_1^{(1)}(x)\right),
\eeq
where $F_2^{(1)}(x)$ is also a solution to
the hypergeometric equation \eqref{HG:1} such that
$$F_2^{(1)}(x)=\ln(1-x)\,F_1^{(1)}(x)+\sum_{n=1}^{\infty}b_n(1-x)^n,$$
with
$$b_n=\frac{(\alpha)_n(\beta)_n}{(n!)^2}\,
\Big( \frac{1}{\alpha}+\cdots+\frac{1}{\alpha+n-1}+
\frac{1}{\beta}+\cdots+\frac{1}{\beta+n-1} - 
2\Big(\frac{1}{1}+\cdots+\frac{1}{n}\Big)\Big).$$
We can check the following $\tau$ is the modulus of the elliptic curve $E_\si$.
\beq\label{mirror-at-lcs}
\tau:=\frac{\pi_{B}(\si)}{\pi_{A}(\si)}.
\eeq
If we put $Q=\exp(2\pi\sqrt{-1}\tau)$, then the $j$-invariant
always has a $Q$-expansion
$$j(\sigma)=\frac{1}{Q}+744+196884Q+\cdots$$
Note that the residue pairing implies
\beq\label{K}
\LD 1,\phi_{-1}\RD=\frac{1}{K(1-C\,\si^l)},
\eeq
where $K$ is some fixed constant. Since the residue pairing
must be identified with the Poincar\'e pairing, the mirror map should satisfy
\beq\label{mm:non-tw}
\Delta_{01}\mapsto 1,\quad \Delta_{02}\mapsto K(1-C\, \si^l) \phi_{-1}(\x)\pi_{A}^2.
\eeq
The next step is to identify the divisor coordinate $t_{02}$ in the
orbifold GW theory and the modulus $\tau$. In order to get the correct $q$-expansion, we define
\beq\label{mm:div}
q:=\exp(t_{-1}), \quad t_{-1}:=\frac{2\pi\sqrt{-1}}{L}\tau,
\eeq
Here $L=3,4,6$ respectively for the elliptic orbifolds $\mathbb{P}^1_{3,3,3},\mathbb{P}^1_{4,4,2},\mathbb{P}^1_{6,3,2}.$ This implies
\beq\label{wronskian}
\LD\Delta_{01},\Delta_{02}\RD\mapsto\LD 1,\frac{\partial}{\partial\,t_{-1}}\RD=\frac{\partial\sigma}{\partial\tau}\frac{\partial\tau}{\partial\,t_{-1}}\LD 1,\frac{\partial}{\partial\sigma}\RD
=\frac{L}{2\pi\sqrt{-1}}\frac{\partial\sigma}{\partial\tau}\frac{1}{\pi_{A}^2}\LD 1,\phi_{-1}\RD
\eeq
is a constant. The last equality follows from equation \eqref{iso:KS} and \eqref{3pt}.  By formulas \eqref{lambda_W}, \eqref{K-constant}, \eqref{mirror-at-lcs}, \eqref{K}, \eqref{mm:div}, we can fix the constant $\lambda_W$ by choosing
$$\LD\Delta_{01},\Delta_{02}\RD\mapsto\LD 1,\frac{\partial}{\partial\,t_{-1}}\RD=1.$$
\begin{remark}
In computations below, for our convenience, we may choose $\pi_{A}(\si)$ different from formula \eqref{lambda_W} by a scalar. This will not change the Frobenius manifold structure since we can always rescale the ring generators to offset its influence.
\end{remark}
\subsubsection{Twisted sectors}
In order to complete the construction, we need to identify the twisted cohomology
classes $\Delta_{ij}$ with monomials $\delta_\r(\si,\x)$. The key observation is that the sections
\beq\label{flat-poly}
\int \delta_\r(\si,\x)\, \frac{\omega}{d W_\si} 
\eeq
of the vanishing cohomology bundle of $W_\si$ are flat with respect to
the Gauss--Manin connection. This way our choice of $\delta_\r$ depends
on an invertible matrix of size $(\mu-2)\times (\mu-2)$. The matrix is decomposed into $1\times1$ and $2\times2$ blocks according to Lemma \ref{PF-equations}.
In particular, the entries of a $2\times2$ block are obtained from a basis of solutions of hypergeometric equations \eqref{2nd-order} with weights $(\alpha_\r,\beta_\r,\gamma_\r)$ near $x=1$:
\beq\label{eq:unity-not0}
\left\{
\begin{aligned}
F_{1,\r}^{(1)}(x) & =  \
_2F_1\left(\alpha_\r,\beta_\r;\alpha_\r+\beta_\r-\gamma_\r+1;1-x\right), \\
F_{2,\r}^{(1)}(x) & = 
\
_2F_1\left(\gamma_\r-\alpha_\r,\gamma_\r-\beta_\r;\gamma_\r-\alpha_\r-\beta_\r+1;1-x\right)\,
(1-x)^{\gamma_\r-\alpha_\r-\beta_\r}.
\end{aligned}
\right.
\eeq

The correlation functions in the
Saito--Givental CohFT are invariant with respect to the translation
$t_{-1}\mapsto t_{-1}+2\pi\sqrt{-1}$, see \eqref{mm:div}. The coefficient in
front of $q^d$, $d\in \Z$, is called the {\em degree-d} part of the
correlator function. By taking the degree-0 part of the 3-point
functions, we obtain a Frobenius algebra structure on the Jacobian algebra
$\mathscr{Q}_W$ that under the mirror map should be identified with
the Frobenius algebra corresponding to the Chen--Ruan orbifold
(classical) cup product.  Using also that
the mirror map preserves homogeneity we obtain a
system of equations for the matrix. It remains only to see that these
equations have a solution. Let us list the explicit formulas for the
mirror map. We omit the details of the
computations, which by the way are best done with the help of some
computer software--Mathematica.

\subsection{\texorpdfstring{Global mirror symmetry for $W_{\sigma}=X_1^3+X_2^3+X_3^3+\sigma X_1X_2X_3$}{LCLS+F}}
The $j$-invariant 
$$j(\si)=-\frac{\si^3(-216+\si^3)^3}{(27+\si^3)^3}=\frac{-27x(8+x)^3}{(1-x)^3}, \quad x=-\frac{\si^3}{27}.$$
$\Phi_{\r}$ satisfies a first order differential equation \eqref{1st-order} and can be solved
$$
\Phi_{\r}=(1-x)^{-\deg\phi_{\r}}\delta_{\r}.
$$
Here $\delta_{\r}$ are some flat sections. The residue pairing implies
\beq\label{residue-E6}
\LD X_1X_2X_3\RD:={\rm Res}\frac{X_1X_2X_3}{(\d_{X_1} W_\si)\, (\d_{X_2} W_\si)\, (\d_{X_3} W_\si)}d^3\x=\frac{1}{27(1-x)}.
\eeq

\subsubsection{GW-point at $x=1$}
For $p_1=-3$, since the weights of equation \eqref{HG:1} in this case are $(\alpha,\beta,\gamma)=(1/3,1/3,2/3)$, according to \eqref{lambda_W}, \eqref{K-constant}, we choose 
$$\pi_{A}(\si)=\ _2F_1\left(\frac{1}{3},\frac{1}{3};1;1-x\right), \quad \pi_B(\si)=\frac{3}{2\pi\sqrt{-1}}\left(F_2^{(1)}(x)-\frac{\ln\, P(p_k)}{3}\cdot F_1^{(1)}(x)\right).$$
According to \eqref{mirror-at-lcs} and \eqref{mm:div}, the Fourier series of $1-x$ in $q=e^{2\pi i\tau/3}$ is
\beq\label{E6-q-expansion-unity}
1-x=-27q-324q^2-2430q^3-13716q^4+\cdots
\eeq 
A natural basis for the flat sections with non-integral degrees are solving from \eqref{1st-order},
$$
\delta_{\r}=(1-x)^{1/3}\phi_{\r}(\x)\, \pi_A,\quad
\delta_{\r'}=(1-x)^{2/3}\phi_{\r'}(\x)\,\pi_A.
$$
Here $\r=(1,0,0),(0,1,0),(0,0,1)$ and $\r'=(1,1,1)-\r.$
Applying \eqref{3pt}, we know the non-vanishing correlators $\LD\cdots\RD_{0,3,0}$ are
$$\LD1,\delta_{\r},\delta_{\r'}\RD_{0,3,0}=\LD\delta_{\r},\delta_{\r},\delta_{\r}\RD_{0,3,0}=\LD\delta_{1,0,0},\delta_{0,1,0},\delta_{0,0,1}\RD_{0,3,0}=\frac{1}{27}.$$ 
The mirror map is given from \eqref{mm:non-tw}, \eqref{mm:div}, with the ring generators identifies by
\beq\label{fermat3-at-unity}
\left(
\begin{matrix}\Delta_{11}\\ \Delta_{21} \\ \Delta_{31}\end{matrix}
\right)\mapsto
\left(
\begin{matrix}1&1&1\\1&\be[\frac{2}{3}]&\be[\frac{1}{3}]\\1&\be[\frac{1}{3}]&\be[\frac{2}{3}]\\
\end{matrix}\right)
\left(
\begin{matrix}\delta_{1,0,0}\\\delta_{0,1,0}\\\delta_{0,0,1}\end{matrix}\right).
\eeq
It is easy to check that this identification agrees with the Chen-Ruan orbifold
cohomology ring of $\mathbb{P}^1_{3,3,3}$ (see 
\eqref{CR-pair}, \eqref{CR-product}). For example, from \eqref{fermat3-at-unity} we have
$$\LD \Delta_{11}, \Delta_{11},\Delta_{11}\RD_{0,3,0}
\mapsto\sum_{\r;\deg\phi_\r=1/3}\LD\delta_{\r},\delta_{\r},\delta_{\r}\RD_{0,3,0}+6\LD\delta_{1,0,0},\delta_{0,1,0},\delta_{0,0,1}\RD_{0,3,0}=\frac{1}{3}.$$
Finally, $\LD\Delta_{11}, \Delta_{21}, \Delta_{31}\RD_{0,3,1}\mapsto 1$ for Lemma \ref{thm:GW recons}
follows from 
\begin{align*}
\LD \Delta_{11}, \Delta_{21}, \Delta_{31}\RD_{0,3}
\mapsto
\frac{(1-(1-x))^{1/3}-1}{9}\pi_{A}=q+q^4+2q^7+\cdots=\frac{\eta(9\tau)^3}{\eta(3\tau)}.
\end{align*}

\subsubsection{GW-point at infinity}
For this limit, the mirror symmetry is already verified in \cite{MR}. 
It maps the twisted sector $\Delta_{i,j}$ to $\delta_{\r}$, where $i$-th index of $\r$ is $j$.

\subsection{\texorpdfstring{Global mirror symmetry for $W_{\sigma}=X_1^4+X_2^4+X_3^2+\sigma X_1^2X_2^2$}{LCLS+7}} The $j$-invariant 
$$
j(\si)=\frac{16(12+\si^2)^3}{(4-\si^2)^2}=\frac{64(3+x)^3}{(1-x)^2}, \quad x=\frac{\si^2}{4}.
$$
For $\r=(10,01,11,21,12)$, $\Phi_{\r}$ satisfies a first
order differential equations \eqref{1st-order}.
On the other hand, both $\Phi_{20}$ and $\Phi_{02}$ satisfy a second order
hypergeometric equation with weights
$\alpha_{20}=\beta_{02}=3/4,\beta_{20}=\alpha_{02}=1/4,\gamma_{20}=\gamma_{02}=1/2$,
which comes from the following system of first order differential equations:
\beq\label{de:2nd-order}
\left\{
\begin{aligned}
(4-\si^2)\, \d_\si\Phi_{20}(\si) & =  \frac{\si}{2}\,
\Phi_{20}(\si)-\Phi_{02}(\si);\\
(4-\si^2)\, \d_\si\Phi_{02}(\si) & =  \frac{\si}{2}\,
\Phi_{02}(\si)-\Phi_{20}(\si).
\end{aligned}
\right.
\eeq
It follows that $\Phi_{02}= L\, \Phi_{20}$ (and $\Phi_{20}= L\,
\Phi_{02})$ where $L$ is the differential operator 
$$L=-(4-\si^2)\d_{\si}+\frac{\si}{2}.$$
Moreover, here the residue pairing implies
\beq\label{residue-E7}
\LD X_1^2X_2^2\RD=\frac{1}{32(1-x)}, \quad 
\LD X_1^4\RD=\LD X_2^4\RD=\frac{-x^{1/2}}{32(1-x)}, \quad 
\LD X_1^3X_2\RD=
\LD X_1X_2^3\RD=0.
\eeq

\subsubsection{GW-point at $x=1$}
We are looking for $\si=2x^{1/2}=\pm2$. We pick $p_2=2.$
The Fourier series for $1-x$ in terms of $q=e^{2\pi i\tau/4}$ is
$$
1-x=64 q^2 - 1536 q^4 + 19200 q^6+\cdots\ .
$$
The weights of \eqref{HG:1} in this case are $(\alpha,\beta,\gamma)=(1/4,1/4,1/2)$. We choose
$$
\pi_{A}(\si)=\ _2F_1\left(\frac{1}{4},\frac{1}{4};1;1-x\right), \quad
\pi_{B}(\si)=\frac{2}{2\pi\sqrt{-1}}\left(F_2^{(1)}(x)-\frac{\ln\, P(p_k)}{2}\cdot F_1^{(1)}(x)\right).
$$
Let us construct a basis of flat sections. For first order equations, 
\beq\label{442-first-order}
\delta_{\r}(\si,\x)=\left(x-1\right)^{\deg\phi_\r}\phi_{\r}(\x)\, \pi_A(\sigma), \quad \r=(10,01,11,21,12)
\eeq
Then for second order equations, we can obtain a pair of polynomials $\delta_{20}$ and $\delta_{02}$
that determine flat sections \eqref{flat-poly} by solving the
following system:
\beq\label{442-second-order}
\left(
\begin{matrix} \phi_{20} \, \pi_A\\ \phi_{02}\, \pi_A\end{matrix} \right)= 
\left(
\begin{matrix}  
\ _2F_1\left(\frac{1}{4},\frac{3}{4};\frac{3}{2};1-x\right) &  
(x-1)^{-1/2}\ _2F_1\left(\frac{1}{4},-\frac{1}{4};\frac{1}{2};1-x\right) \\
\ _2F_1\left(\frac{1}{4},\frac{3}{4};\frac{3}{2};1-x\right)&  
 -(x-1)^{-1/2}\ _2F_1\left(\frac{1}{4},-\frac{1}{4};\frac{1}{2};1-x\right)
\end{matrix}
\right) 
\left( \begin{matrix} \delta_{20} \\ \delta_{02}\end{matrix} 
\right).
\eeq
Now the genus-0 3-point Saito-Givental correlators for flat sections $\delta_{\r}$ can be calculated using residue formula \eqref{residue-E7}. Based on this, we can check that after a linear transformation, we actually match those flat elements $\Delta_{i,j}$ in Chen-Ruan cohomology of $\mathbb{P}_{4,4,2}$ to the flat sections via the following mirror map 
\ben\label{442-mirror-unity}
\left\{
\begin{array}{lll}
\Delta_{1,1}\mapsto\delta_{10}+\sqrt{-1}\delta_{01}, &\Delta_{1,2}\mapsto 2\sqrt{-1}\delta_{11}+2\delta_{02}, &\Delta_{1,3}\mapsto 4\sqrt{-1}\delta_{2,1}-4\delta_{1,2};\\
\Delta_{2,1}\mapsto\sqrt{-1}\delta_{10}+\delta_{01}, &\Delta_{2,2}\mapsto 2\sqrt{-1}\delta_{11}-2\delta_{02}, &\Delta_{2,3}\mapsto -4\delta_{2,1}+4\sqrt{-1}\delta_{1,2};\\
\Delta_{01}\mapsto 1, & \Delta_{3,1}\mapsto 8\delta_{20}, & \Delta_{02}\mapsto 32(1-x)\, \phi_{-1}\,\pi_A^2.
\end{array}
\right.
\een
It is not hard to check this map matches all the pairing \eqref{CR-pair} using \eqref{residue-E7}. For example, 
$$\LD 8\delta_{20},8\delta_{20}\RD=16\LD\phi_{20}+\phi_{02},\phi_{20}+\phi_{02}\RD\ \ _2F_1\left(\frac{1}{4},\frac{3}{4};\frac{3}{2};1-x\right)^{-2}=\frac{1}{2}.$$
The last equality is a consequence of the hypergeometric identity
$$\ _2F_1\left(\frac{1}{4},\frac{3}{4};\frac{3}{2};1-x\right)^2(1+x^{1/2})=2.$$
We can also check the mirror map matches Chen-Ruan product \eqref{CR-product}. For example,
$$\begin{aligned}
\LD \Delta_{1,1}, \Delta_{1,1}, \Delta_{1,2}\RD_{0,3}
\mapsto&\frac{1}{8}\left(\ _2F_1\left(\frac{1}{4},-\frac{1}{4};\frac{1}{2};1-x\right)+1\right)\ _2F_1\left(\frac{1}{4},\frac{1}{4};1;1-x\right)\\
&=\frac{1}{4}+q^4+q^{8}+q^{16}+\cdots.
\end{aligned}
$$
The last equality follows from mirror map \eqref{mirror-at-lcs}  and a hypergeometric identity 
$$\ _2F_1\left(\frac{1}{4},-\frac{1}{4};\frac{1}{2};1-x\right)^2=\frac{1+x^{1/2}}{2}.$$
The degree-1 genus-0 3-point correlator in Lemma \ref{thm:GW recons} is verified by
\ben
\begin{aligned}
\LD \Delta_{11}, \Delta_{21}, \Delta_{31}\RD_{0,3}
\mapsto&\frac{1}{8}\left(1-x\right)^{1/2}\ _2F_1\left(\frac{1}{4},\frac{3}{4};\frac{3}{2};1-x\right)\ _2F_1\left(\frac{1}{4},\frac{1}{4};1;1-x\right)\\
&=q+2q^5+q^{9}+2q^{13}+\cdots.
\end{aligned}
\een

\subsubsection{GW-point at infinity}
Near the point $\sigma=\infty$, we choose
$$
\pi_{A}(\si)=\left(\frac{1}{4x}\right)^{1/4}\ _2F_1\left(\frac{1}{4},\frac{3}{4};1;\frac{1}{x}\right), \quad
\pi_B(\si)=\frac{2}{2\pi\sqrt{-1}}\left(F_2^{(1)}(x)-\frac{\ln\, P(p_k)}{2}\cdot F_1^{(1)}(x)\right),
$$
$\pi_{A}(\si)$ and $\pi_{A}(\si)$ satisfy \eqref{HG:1}. $F_2^{(1)}(x)$ is defined by (same $b_n$ as in \eqref{K-constant}), 
$$F_2^{(1)}(x)=-(\ln x)\,F_1^{(1)}(x)+\sum_{n=1}^{\infty}b_n\,x^{-n}.$$
We construct a mirror map \eqref{mirror-at-lcs} for the K\"ahler class. It implies an Fourier expansion
$$x^{-1}=64q^4-2560q^8+84736q^{12}+...$$  
We choose the flat sections for the first order differential equations
$$\delta_{\r}=(x-1)^{\deg\phi_{\r}}\phi_{\r}\,\pi_{A}.$$
For the second order differential equations, we solve \eqref{de:2nd-order} to obtain
\beq
\left(
\begin{matrix} \phi_{20} \, \pi_A\\ \phi_{02}\, \pi_A\end{matrix} \right)= 
\left(
\begin{matrix} 
 x^{-3/4}\ _2F_1\left(\frac{3}{4},\frac{5}{4};\frac{1}{2};\frac{1}{x}\right)&
x^{-1/4}\ _2F_1\left(\frac{1}{4},\frac{3}{4};\frac{1}{2};\frac{1}{x}\right)\\  
 -2x^{-1/4}\ _2F_1\left(\frac{1}{4},\frac{3}{4};\frac{1}{2};\frac{1}{x}\right) &
 -\frac{1}{2}x^{-3/4}\ _2F_1\left(\frac{3}{4},\frac{5}{4};\frac{1}{2};\frac{1}{x}\right)
\end{matrix}
\right) 
\left( \begin{matrix} \delta_{20} \\ \delta_{02}\end{matrix} 
\right).
\eeq
The entries comes from the solutions of a hypergeometric equation \eqref{2nd-order} near $x=\infty$:
\beq\label{eq:infinity-not0}
\left\{
\begin{aligned}
F_{1,\r}^{(\infty)}(x)& = 
\
_2F_1\left(\alpha_\r,\alpha_\r-\gamma_\r+1;\alpha_\r-\beta_\r+1;x^{-1}\right)\,
x^{-\alpha_\r} \, ,\\
F_{2,\r}^{(\infty)}(x) & = 
\
_2F_1\left(\beta_\r,\beta_\r-\gamma_\r+1;\beta_\r-\alpha_\r+1;x^{-1}\right)\,
x^{-\beta_\r}\, .
\end{aligned}
\right.
\eeq

We construct a mirror map as follows
\ben\label{E7-mirror-infinity}
\left\{
\begin{array}{lll}
\Delta_{1,1}\mapsto2\delta_{10}, &\Delta_{1,2}\mapsto-4\sqrt{2}\delta_{02}, &\Delta_{1,3}\mapsto-8\delta_{1,2};\\
\Delta_{2,1}\mapsto2\sqrt{-1}\delta_{01}, &\Delta_{2,2}\mapsto-2\sqrt{2}\delta_{20}, &\Delta_{2,3}\mapsto8\sqrt{-1}\delta_{2,1};\\
\Delta_{01}\mapsto1, & \Delta_{3,1}\mapsto4\sqrt{-1}\delta_{11}, &\Delta_{02}\mapsto32(1-x)\, X_1^2X_2^2\,\pi_A^2.
\end{array}
\right.
\een
Now, everything can be checked as before using Mathematica. In particular, we get 
$$
\LD \Delta_{1,1}, \Delta_{2,1}, \Delta_{3,1}\RD_{0,3}
\mapsto \frac{\pi_{A}(\si)}{2}
=q+2q^{5}+q^{9}+\cdots.
$$

\subsection{\texorpdfstring{Global mirror symmetry for $W_{\sigma}=X_1^6+X_2^3+X_3^2+\sigma X_1^4X_2$}{LCLS+E8}}
The $j$-invariant 
\ben
j(\si)=1728\frac{4\si^3}{27+4\si^3}=\frac{-1728x}{1-x}, \quad x=-\frac{4}{27}\si^3.
\een
In order to construct the mirror map for twisted sectors (see Table \ref{weight-E8}), we have to find
a basis of homogeneous flat sections (with non-integer degrees) of the
Gauss--Manin connection. $\Phi_{11}$ and $\Phi_{15}$ satisfy first order equations and
the periods corresponding to the polynomials $\phi_{(k,1)}$ and $\phi_{(k+2,0)}$ satisfy
\beq\label{ode-E8}
\left\{
\begin{aligned}
(27+4\si^3)\d_\si \, \Phi_{(k+2,0)} & = -(k+3)\si^2\, \Phi_{(k+2,0)} -\frac{9(k+1)}{2}\Phi_{(k,1)} \\
(27+4\si^3)\d_\si \, \Phi_{(k,1)} & =  \frac{3(k+3)\si}{2}\Phi_{(k+2,0)} -  (k+1)\si^2 \Phi_{(k,1)} 
\end{aligned}
\right.
\eeq
Moreover, we know $\Phi_{(k,1)}$
satisfies a hypergeometric equation (see Section \ref{PF:system}). Let 
\beq\label{diff-op-E8}
L_{k}:=\frac{2}{3(k+3)\si}\left( (27+4\si^4)\d_{\si}+(k+1)\si^2\right),\quad k=0,1,2
\eeq 

\subsubsection{GW-point at $x=1$}
Near $p_2=\frac{3}{2}(-2)^{1/3},$ $N=1$. The weights of \eqref{HG:1} in this case are $(\alpha,\beta,\gamma)=(1/12,7/12,2/3)$. We choose
$$
\pi_{A}(\si)=\ _2F_1\left(\frac{1}{12},\frac{7}{12};1;1-x\right), \quad
\pi_{B}(\si)=\frac{1}{2\pi\sqrt{-1}}\left(F_2^{(1)}(x)-\ln\, P(p_k)\cdot F_1^{(1)}(x)\right).
$$
The mirror map for K\"ahler class \eqref{mirror-at-lcs} implies Fourier expansion
$$1-x=-1728q^6-1700352q^{12}+\cdots$$ 
We choose the following map from the non-twisted sectors in GW theory to the flat sections \eqref{mm:non-tw}, 
\ben
\Delta_{0,1}\mapsto 1, \quad \Delta_{02}\mapsto36(1-x)X_1^4X_2\pi_A^2.
\een
For the twisted sectors, we choose flat sections up to some constant (denoted by $\sim$). For first order differential equations \eqref{1st-order}, according to Table \ref{weight-E8}, we choose
$$\Delta_{11}\sim (1-x)^{1/6}\, \phi_{10} \, \pi_A, \quad \Delta_{15}\sim (1-x)^{5/6}\, \phi_{31}\, \pi_A.$$
Others are obtained by solving equations \eqref{ode-E8} (or \eqref{2nd-order}) near $x=1$,
\ben
\left\{
\begin{aligned}
&
\Delta_{21}\sim (1-x)^{1/3}\Big(F_{2,20}^{(1)}(x)\, \phi_{01}+(-2)^{-1/3}\
F_{2,01}^{(1)}(x)\, \phi_{20}\Big) \, \pi_A,\\
& \Delta_{12}\sim (1-x)^{1/3}\Big(F_{1,20}^{(1)}(x)\, \phi_{01}-3(-2)^{-1/3}\ F_{1,01}^{(1)}(x)\, \phi_{20}\Big) \, \pi_A,\\
& \Delta_{31}\sim (1-x)^{1/2}\Big(F_{2,30}^{(1)}(x)\, \phi_{11}+(-2)^{-1/3}\ F_{2,11}^{(1)}(x)\, \phi_{30}\Big) \, \pi_A,\\
& \Delta_{13}\sim (1-x)^{1/2}\Big(F_{1,20}^{(1)}(x)\, \phi_{11}-2(-2)^{-1/3}\ F_{1,11}^{(1)}(x)\, \phi_{30}\Big) \, \pi_A,\\
& \Delta_{22}\sim (1-x)^{2/3}\Big(F_{2,40}^{(1)}(x)\, \phi_{21}+(-2)^{-1/3}\ F_{2,21}^{(1)}(x)\, \phi_{40}\Big) \, \pi_A,\\
& \Delta_{14}\sim (1-x)^{2/3}\Big(F_{1,40}^{(1)}(x)\, \phi_{21}-\frac{5}{3}(-2)^{-1/3}\ F_{1,21}^{(1)}(x)\, \phi_{40}\Big) \, \pi_A.\\
\end{aligned}\right.
\een 
$F_{i,\r}^{(1)}(x)$ is from \eqref{eq:unity-not0} with weights in Table \ref{weight-E8}.
The propotions can be fixed by identify the pairng and the ring structure constants. For the 3-point correlator we get 
\begin{eqnarray*}
\LD \Delta_{1,1}, \Delta_{2,1}, \Delta_{3,1}\RD_{0,3}
\mapsto q+2q^7+\cdots
\end{eqnarray*}

\subsubsection{FJRW-point at infinity} 

The hypergeometric equation for the periods $\pi_{A}$ has weights
$(\alpha,\beta,\gamma)=(1/12,7/12,2/3)$. 
Since $\alpha-\beta$ is not an integer, the monodromy is
diagonalizable and we have the following basis of solutions
(eigenvectors for the monodromy around $\si=\infty$)  near
$x=\infty:$
\ben
\begin{aligned}
\pi_{A_{\infty}}:=x^{-1/12}\
_2F_1\Big(\frac{1}{12},\frac{5}{12};\frac{1}{2};x^{-1}\Big),\quad 
\pi_{B_{\infty}}:=\la\ x^{-7/12}\ _2F_1\Big(\frac{7}{12},\frac{11}{12};\frac{3}{2};x^{-1}\Big),
\end{aligned}
\een
where the constant $\la$ will be fixed later on. 
Put 
\beqa
t_{-1}=\pi_{B_{\infty}}/\pi_{A_{\infty}} \approx \la\ x^{-\frac{1}{2}},
\eeqa
where $\approx$ means that we truncated terms of order $O(\si^2)$. 
It is easy to check (by using the differential equation for the periods) that when restricted to the subspace of marginal
deformation, $t_{-1}$ is a degree $0$ flat coordinate, i.e., the residue
pairing $\langle 1, \d/\d t_{-1}\rangle$ is a constant. 

For first order differential equations \eqref{1st-order}, from Table \ref{weight-E8}, we get flat sections
$$A_{(1,0)}=(x-1)^{1/6}\phi_{10}\pi_{A}(\si), \quad A_{(3,1)}=(x-1)^{5/6}\phi_{31}\pi_{A}(\si).$$ 
The solutions to the differential equations \eqref{ode-E8} near $x=\infty$, are
\ben
\left(
\begin{matrix}
\phi_{(k,1)}\,\pi_A(\si)\\
\phi_{(k+2,0)}\,\pi_A(\si)
\end{matrix}
\right)
=
\left(
\begin{matrix}
F_{1,(k,1)}^{(\infty)}(x)& F_{2,(k,1)}^{(\infty)}(x)\\
L_{k} F_{1,(k,1)}^{(\infty)}(x)&L_{k} F_{2,(k,1)}^{(\infty)}(x)
\end{matrix}\right)
\left(
\begin{matrix}
A_{(k,1)}\\
A_{(k+2,0)}
\end{matrix}
\right),\quad k=0,1,2
\een
Here $F_{i,\r}^{(\infty)}(x)$ is from \eqref{eq:infinity-not0} and $L_k$ is from \eqref{diff-op-E8}.
Let $c_{\r}$ be given below
\begin{table}[h]\center
\begin{tabular}{|c|c|c|c|c|c|c|c|c|}
\hline

$\r=$ & $10$&$01$ & $20$ & $11$ & $30$ & $21$ & $40$&$31$

\\
\hline
$c_{\r}$ & $\la_1$  & $\la_2$& $-\frac{\la_1^2C_0}{3}$ & $\la_1\la_2$&$-\la_1^3C_0$ & $\la_1\la_2^2$& $\frac{4\la_1^4C_0}{5}$&$\frac{2\la_1^5C_0}{9}$

\\
\hline
\end{tabular}
\end{table}

The constants appearing in the table are given as follows:
\beq\label{const}
\la_1^6=24C_0^2,\quad \la_2^2=\frac{\la_1^4}{C_0}, \quad C_0^3=-\frac{27}{4}.
\eeq

Let $\delta_\r(\s,\x)$ be polynomials, such that the geometric sections satisfying
(see \eqref{s-geom}) 
$$[\delta_{\r}\omega]=c_{\r}A_{\r}.$$
Now we compute the pairing and the necessary genus-$0$ correlators. 
The pairing is 
$$
\LD\delta_{10},\delta_{31}\RD=\LD\delta_{01},\delta_{21}\RD=\LD\delta_{20},\delta_{40}\RD=\LD\delta_{11},\delta_{11}\RD=\LD\delta_{30},\delta_{30}\RD=1.
$$
All 3-point correlator functions that do not have insertion 1
(otherwise the correlator reduces to a 2-point one) have a limit at
$\si=\infty$. The non-zero limits are as follows:
\ben
\begin{aligned}
&\LD\delta_{10},\delta_{10},\delta_{40}\RD_{0,3}
=\LD\delta_{10},\delta_{01},\delta_{11}\RD_{0,3}
=\LD\delta_{01},\delta_{01},\delta_{20}\RD_{0,3}
=1.\\
&\LD\delta_{10},\delta_{20},\delta_{30}\RD_{0,3}=\LD\delta_{20},\delta_{20},\delta_{20}\RD_{0,3}=-3.
\end{aligned}
\een
In other words, the Jacobian algebra extends over $\si=\infty$. If we
denote the extension by $\mathscr{Q}_{W^{\infty}}$, then it is not
hard to see that $\delta_{10}$ and $\delta_{01}$ are generators and we
have
\ben
\mathscr{Q}_{W^{\infty}}:=\C[\delta_{10},\delta_{01}]/\left(4\delta_{10}^3\delta_{01},\delta_{10}^4+3\delta_{01}^2\right).
\een
Finally, we set
\ben
\la=\frac{\la_1^4\la_2}{54}.
\een 
The nonzero 4-point genus-$0$ basic correlators are
\begin{align*}
\LD\delta_{01},\delta_{01},\delta_{01},\delta_{-1}\RD_{0,4}
&=\frac{\partial}{\partial\,t_{-1}}\LD\la_2^3\left(x^{1/12}\Phi_{01}+\frac{3}{4C_0}x^{-1/4}\Phi_{20}\right)^{3}\RD\Big\vert_{x=\infty}
\\
&=-\frac{\la_2^2C_0}{4\la_1^4}\frac{\d}{\d\,t_{-1}}\left(\la\,x^{-1/2}\right)=-\frac{1}{4}.
\end{align*}
and
$$\LD\delta_{10},\delta_{10},\delta_{10}^2\delta_{01},\delta_{-1}\RD_{0,4}
=-\frac{1}{4}.$$
Recall that $\rho_1,\rho_{2}$ are generators of the FJRW ring
for $\Big(W'=X_1^4X_2+X_2^3+X_3^2,
G_{W'}\Big).$
We construct the following mirror map from $\mathscr{H}_{W'}$ to $\mathscr{Q}_{W^{\infty}}$,
\ben
\Big(\rho_1,\rho_2\Big)\mapsto\Big(
(-1)^{5/6}\delta_{01},(-1)^{2/3}\delta_{10}\Big)
\een
Using Lemma \ref{lm:recons-FJRW}, it is easy to check that this map identifies the FJRW theory of $(W', G_{W'})$ to the Saito-Givental limit of $W_{\si}=X_1^6+X_2^3+X_3^2+\si X_1^4X_2$ at $\si=\infty,$
$$\mathscr{A}_{W'}^{\rm FJRW}=\lim_{\sigma\to\infty}\mathscr{A}_W^{\rm SG}(\si).$$


\begin{thebibliography}{JKV2}


\bibitem{A} Acosta, Pedro:
{\em FJRW-rings and Landau-Ginzburg mirror symmetry in two dimensions.}
arXiv:0906.0970 

\bibitem{AGV}{Arnold, V.; Gusein-Zade, S.; Varchenko, A.:}
{\em Singularities of Differentiable maps.} Vol. II. Monodromy and
  Asymptotics of Integrals. Boston, MA: Birkh\"auser Boston,
  1988. viii+492pp

\bibitem{BCOV}
Bershadsky, M.; Cecotti, S.; Ooguri, H.; Vafa, C.:
{\em Kodaira-Spencer theory of gravity and exact results for quantum string amplitudes.}
Comm. Math. Phys. 165 (1994), no. 2, 311-427

\bibitem{BH} Berglund, Per; H\"ubsch, Tristan:
{\em A generalized construction of mirror manifolds.} Nuclear Phys. B 393 (1993), no. 1-2, 377-391.


\bibitem{CDGP} Candelas, Philip; de la Ossa, Xenia C.; Green, Paul S.; Parkes, Linda:
{\em A pair of Calabi-Yau manifolds as an exactly soluble superconformal theory.} Nuclear Phys. B 359 (1991), no. 1, 21-74.

\bibitem{CheR}{Chen, Weimin; Ruan, Yongbin:}
\emph{ Orbifold Gromov-Witten theory.} Orbifolds in mathematics and physics.
Contemp. Math., 310, Amer. Math. Soc., Providence, RI(2002): 25--85.

\bibitem{C} Chiodo, Alessandro:
{\em Towards an enumerative geometry of the moduli space of twisted curves and $r$-th roots.}
Compos. Math. 144 (2008), no. 6, 1461-1496.

\bibitem{CIR} Chiodo, Alessandro; Iritani, Hiroshi; Ruan, Yongbin:
{\em Landau-Ginzburg/Calabi-Yau correspondence, global mirror symmetry and Orlov equivalence.} 
Preprint arXiv:1201.0813v3. To appear in Publ. Math. Inst. Hautes $\acute{\rm E}$tudes Sci.


\bibitem{CR1} Chiodo, Alessandro; Ruan, Yongbin:
{\em A global mirror symmetry framework for the Landau-Ginzburg/Calabi-Yau correspondence.} 
Ann. Inst. Fourier (Grenoble) vol. 61, no. 7 (2011), 2803-2864.

\bibitem{CR2} Chiodo, Alessandro; Ruan, Yongbin:
{\em Landau-Ginzburg/Calabi-Yau correspondence for quintic three-folds via symplectic transformations.}
Invent. Math. 182 (2010), no. 1, 117-165.

\bibitem{CI} Coates, Tom; Iritani, Hiroshi:
{\em On the Convergence of Gromov-Witten Potentials and Givental's Formula}
Preprint arXiv:1203.4193v1.

\bibitem{CL} 
Costello, Kevin; Li, Si:
{\em Quantum BCOV theory on Calabi-Yau manifolds and the higher genus B-model.}
Preprint arXiv:1201.4501v1

\bibitem{DM}
Doran, Charles F.; Morgan, John W.:
{\em Mirror symmetry and integral variations of Hodge structure underlying one-parameter families of Calabi-Yau threefolds.} 
Mirror symmetry. V, 517-537, AMS/IP Stud. Adv. Math., 38, Amer. Math. Soc., Providence, RI, 2006.

\bibitem{Du}{Dubrovin, Boris:}
\emph{Geometry of 2d Topological Field Theories.} 
Integrable Systems and Quantum Groups. Lecture Notes in Math. 1620: Springer, Berlin(1996): 120-348


\bibitem{FJR} Fan, Huijun; Jarvis, Tyler J.; Ruan, Yongbin:
{\em The Witten equation, mirror symmetry, and quantum singularity theory.} 
Ann. of Math. (2) 178 (2013), no. 1, 1-106.

\bibitem{FJR1} Fan, Huijun; Jarvis, Tyler J.; Ruan, Yongbin:
{\em The Witten Equation and Its Virtual Fundamental Cycle.} 
Book in preparation. Preprint arXiv:0712.4025v3.

\bibitem{FS} Fan, Huijun; Shen, Yefeng:
{\em Quantum ring of singularity $X^{p}+XY^{q}$}. 
Michigan Math. J. 62 (2013), no. 1, 185-207. 


\bibitem{G}{G\"ahrs, Swantje}.
\emph{Picard-Fuchs equations of special one-parameter families of invertible polynomials.} 
Fields Institute Communications, Volume 67, 2013, 285-310.

\bibitem{G0} 
Givental, Alexander B.:
{\em Equivariant Gromov-Witten invariants.}
Internat. Math. Res. Notices 1996, no. 13, 613-663.

\bibitem{G1} Givental, Alexander B.:
{\em Semisimple Frobenius structures at higher genus.}
Internat. Math. Res. Notices 2001, no. 23, 1265-1286.

\bibitem{G2} Givental, Alexander B.:
{\em Gromov-Witten invariants and quantization of quadratic Hamiltonians.} Dedicated to the memory of I. G. Petrovskii on the occasion of his 100th anniversary. Mosc. Math. J. 1 (2001), no. 4, 551-568, 645. 


\bibitem{He} Hertling, Claus:
\emph{Frobenius Manifolds and Moduli Spaces for Singularities.}
Cambridge Tracts in Mathematics, 151. Cambridge University Press, Cambridge, 2002. x+270 pp. 

\bibitem{HKQ}
Huang, M.-x.; Klemm, A.; Quackenbush, S.:
{\em Topological string theory on compact Calabi-Yau: modularity and boundary conditions.} Homological mirror symmetry, 45-102, Lecture Notes in Phys., 757, Springer, Berlin, 2009.


\bibitem{K}  Krawitz, Marc:
{\em FJRW rings and Landau-Ginzburg Mirror Symmetry.}
Ph.D. thesis, University of Michigan, Ann Arbor, MI, 2010.

\bibitem{KS}{Krawitz, Marc; Shen, Yefeng}:
{\em  Landau-Ginzburg/Calabi-Yau correspondence of all genera for elliptic orbifold $\mathbb{P}^1$.}
Preprint arXiv:1106.6270

\bibitem{L}{Li, Si}: {\em On the quantum theory of Landau-Ginzburg B-model.} Preprint.

\bibitem{LLY}
Lian, Bong H.; Liu, Kefeng; Yau, Shing-Tung:
{\em Mirror principle. I.} Asian J. Math. 1 (1997), no. 4, 729-763.

\bibitem{LMS} {Liu, Sijun; Milanov, Todor; Shen, Yefeng}:
{\em Global mirror symmetry for invertible simple elliptic singularities (II)}.
In preparation.



\bibitem{M}{Milanov, Todor}: 
{\em Analyticity of the total ancestor potential in singularity theory}.
Preprint arXiv:1303.3103

\bibitem{MR}{Milanov, Todor; Ruan, Yongbin}:
\emph{Gromov-Witten theory of elliptic orbifold $\mathbb{P}^1$ and quasi-modular forms.}
Preprint arXiv:1106.2321v1.

\bibitem{MRS}{Milanov, Todor; Ruan, Yongbin; Shen, Yefeng}:
\emph{
Gromov--Witten theory and cycled-valued modular forms.}
Preprint arXiv:1206.3879v1.

\bibitem{MS}
{Milanov, Todor; Shen, Yefeng}:
\emph{The modular group for the total ancestor potential of a simple elliptic singularity.}
Preprint.


\bibitem{Ru1} Ruan, Yongbin:
{\em The Witten equation and the geometry of the Landau-Ginzburg model.} String-Math 2011, 209-240, Proc. Sympos. Pure Math., 85, Amer. Math. Soc., Providence, RI, 2012.


\bibitem{S0} Saito, Kyoji:
{\em Primitive forms for a universal unfolding of a function with an isolated critical point.} 
J. Fac. Sci. Univ. Tokyo Sect. IA Math. 28 (1981), no. 3, 775-792 (1982).


\bibitem{S2} Saito, Kyoji:
\emph{Einfach-elliptische Singularit\"aten}. 
Invent. Math. 23(1974): 289--325


\bibitem{SaT} Saito, Kyoji; Takahashi, Atsushi:
{\em From primitive forms to Frobenius manifolds.}
From Hodge theory to integrability and TQFT tt*-geometry, 31-48, Proc. Sympos. Pure Math., 78, Amer. Math. Soc., Providence, RI, 2008.

\bibitem{ST} Satake, Ikuo; Takahashi, Atsushi:
{\em Gromov–Witten invariants for mirror orbifolds of simple elliptic singularities.}
Ann. Inst. Fourier (Grenoble) vol. 61, no. 7 (2011), 2885-2907.

\bibitem{S} Shen, Yefeng:
{\em Gromov-Witten theory of elliptic orbifold projective lines.}
Ph.D. thesis, University of Michigan, Ann Arbor, MI, 2013.

\bibitem{Te}
Teleman, Constatin:
\emph{The structure of 2D semisimple field theories.}
Invent. Math. Volume 188, Number 3 (2012), 525-588.


\end{thebibliography}
\end{document}